# Quantales and Hyperstructures

Monads, Mo' Problems

by

Andrew Joseph Dudzik

A dissertation submitted in partial satisfaction of the

requirements for the degree of

Doctor of Philosophy

in

Mathematics

in the

Graduate Division

of the

University of California, Berkeley

Committee in charge:

Professor Martin Olsson (Chair)
Professor Matt Baker (Advisor)
Professor Allan Sly

Fall 2016

Abstract

Quantales and Hyperstructures

by

Andrew Joseph Dudzik

Doctor of Philosophy in Mathematics

University of California, Berkeley

Professor Martin Olsson (Chair)

We present a theory of lattice-enriched semirings, called *quantic semirings*, which generalize both quantales and powersets of hyperrings. Using these structures, we show how to recover the spectrum of a Krasner hyperring (and in particular, a commutative ring with unity) via universal constructions, and generalize the spectrum to a new class of hyperstructures, *hypersemirings*. (These include hyperstructures currently studied under the name "semihyperrings", but we have weakened the distributivity axioms.)

Much of the work consists of background material on closure systems, suplattices, quantales, and hyperoperations, some of which is new. In particular, we define the category of covered semigroups, show their close relationship with quantales, and construct their spectra by exploiting the construction of a universal quotient frame by Rosenthal.

We extend these results to hypersemigroups, demonstrating various folkloric correspondences between hyperstructures and lattice-enriched structures on the powerset. Building on this, we proceed to define quantic semirings, and show that they are the lattice-enriched counterparts of hypersemirings. To a quantic semiring, we show how to define a universal quotient quantale, which we call the *quantic spectrum*, and using this, we show how to obtain the spectrum of a hypersemiring as a topological space in a canonical fashion.

Finally, we we conclude with some applications of the theory to the ordered blueprints of Lorscheid.

Wer das Tiefste gedacht, liebt das Lebendigste.

— Friedrich Hölderlin,
*Sokrates und Alcibiadese*

In memory of Kiel Sturm

1983 – 2012

## ACKNOWLEDGMENTS

I could not possibly have completed this work without the support of many amazing people, who are too numerous to list without fear of omission. But thank you Dad, and thank you Anne—you were there for the worst of it. And thanks to the old friends who believed in me for no apparent reason.

Thank you Ellen; you sent pie. Broad thanks to the rest of my family, to my fellow students, and to the co-workers who respected my decision to go back to school, even while warning me how hard it would be.

I would have not known where to begin without my advisor Matt Baker, who repeatedly challenged me to produce better work than I thought possible. His patience with my struggles and false starts was more valuable than any mathematical contribution.

Oliver Lorscheid and Ben Steinberg provided much-needed feedback on an early draft of this work. Oliver helped give me confidence that the language of closure operators is a good way to study exotic algebraic structures, and Ben helped give me confidence that applying quantales to the theory of hyperrings is not completely crazy.

Finally, thanks to the person who opened the doors for me to pursue mathematics in the first place, Zvezdelina Stankova. When I was in high school, her work at the Berkeley Math Circle provided me with the backbone for all of my future explorations of mathematics. When I was in graduate school, working as her teaching assistant was a master class in running a successful university course.

CONTENTS

# INTRODUCTION 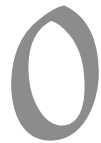

> We cannot understand what something is without grasping what, under certain conditions, it can become.
>
> — Roberto Mangabeira Unger,
> *The Singular Universe and the Reality of Time*

Classical algebra, in its most innocent form, concerns the study of sets, functions, and operations $A^n \to A$ of finite arity $n$, which are usually required to satisfy some list of axioms. In this generality, we encounter all the familiar objects, like groups, rings, and fields, as well as their less mainstream counterparts, including semigroups, semirings, and so on.

There is a growing fashion to consider a more general kind of algebra, *hyperalgebra*, in which some functions $A \to B$ (usually the operations) are replaced by *hyperfunctions* $A \to \mathcal{P}(B)$ taking values in the powerset of $B$. Classical algebra then becomes the special case of hyperalgebra in which all morphisms and all operations take only singleton values.

There are tremendous *pro tanto* benefits to the study of hyperalgebra, by which we mean that there are very real advantages that must be weighed against the increased difficulty of a more general theory. We will begin by listing some of these advantages, and their manifestations in the developing theory of *Krasner hyperrings*—a hyper-analogue of commutative rings with unity which still have single-valued multiplication, but which are allowed to have multi-valued addition.

The first advantage is one of unified language, of which tropical bend relations are a key example. If $p(\vec{X})$ is a polynomial with coefficients in the tropical semifield $\mathbb{T} = (\mathbb{R}_{\geqslant 0}, \cdot, \max)$, then the tropical vanishing locus of $p$ has a definition that is not always easy to work with: it is the set of points $\vec{a}$ where the minimum value of the monomials in $p(\vec{a})$ is achieved at least twice. However, if we instead regard $\mathbb{T}$ as a *hyperfield*—a field where the addition is allowed to be multi-valued—by setting $a + b = \{\max(a, b)\}$ for $a \neq b$ and $a + a = [0, a]$, then $\vec{a}$ is in the tropical vanishing locus if and only if $0 \in p(\vec{a})$.[17] In this setting and others, hyperfields allow us to place new concepts and old ones beneath the same umbrella.

The second advantage is the formation of arbitrary quotients. If $A$ is a set equipped with an operation $* : A^n \to A$, then identifying elements of $A$ according to an arbitrary equivalence relation will generally not allow us to define an operation on the quotient. However, if $*$ is a hyperoperation, there is an inherited hyperoperation: if $A_1, \cdots, A_n$ are equivalence classes, just define $*(A_1, \cdots, A_n)$ to be the set of equivalence classes generated by $\bigcup_{a_i \in A_i} *(a_1, \cdots, a_n)$. This is the perspective of Connes and Consani in studying the hyperring of adèle classes [2]—the adèles form a ring, and the adèle classes form a multiplicative group, but the addition becomes a hyperaddition in the quotient.

The final advantage we'll mention has to do with the surprising connections between hyperfields and projective geometry, which hint at the broad class of phenomena lurking within hyperstructures. For example, there is a close correspondence between projective incidence groups and extensions of the *Krasner hyperfield* [4], which is the final object in the category of Krasner hyperrings. In the ongoing related work of Baker and Bowler [1], a definition of matroids over hyperfields is used to unify seemingly disparate theories of matroids, such as matroids over the tropical hyperfield being precisely valuated matroids.

As promised, these advantages come with technical challenges relative to classical algebra. Many of these challenges arise from one fundamental problem: sets are basically discrete, while powersets have structure, including a natural partial order. In practice, this means that axioms and propositions will tend to involve the asymmetric predicate $\subset$ rather than the symmetric predicate $=$. For example, we mentioned above that the condition $\{0\} \subset f(\vec{a})$ is the appropriate hyper-analogue of vanishing of polynomials. But the fact that we cannot usually write the opposite inequality $f(\vec{a}) \subset \{0\}$ introduces severe difficulties in adapting classical proofs to the hyper-setting.

Another mathematical mainstay threatened by this asymmetry is the *distributive law*, which does not lift to the powerset. Instead, we often have to make do with the *subdistributive law*:

$$S \cdot (T + U) \subset S \cdot T + S \cdot U$$

Consider the following example in $\mathcal{P}(\mathbb{Z})$:

$$\{0, 2\} = \{0, 1\} \cdot (\{1\} + \{1\})$$
$$\subsetneq \{0, 1\} \cdot \{1\} + \{0, 1\} \cdot \{1\} = \{0, 1, 2\}$$

Another naturally-occurring example of sub-distributivity is the arithmetic of intervals, an area of research that has acknowledged sub-distributivity for at least a half-century, e.g.:

$$[-1, 1] = [0, 1] \cdot ([-2, -1] + [1, 2])$$
$$\subsetneq [0, 1] \cdot [-2, -1] + [0, 1] \cdot [1, 2] = [-2, 2]$$

The loss of distributivity is particularly harmful for the theory of polynomials. If a hyperring A is not distributive, then the collection of polynomials A[X], with its natural hyperoperations, will fail to have associative multiplication.

A challenge of a more philosophical nature arises when we attempt to generalize basic concepts from ring theory, such as ideals and prime spectra. Jun [7] has given a description of the spectrum of prime ideals of a Krasner hyperring, and showed that—at least in the integral case—it shares many of the nice properties of its classical counterpart. However, in generalizing this construction further (say, to multi-valued multiplication), it is difficult to know in advance whether we should expect our imitations of the classical definitions to hold up. In practice, simply porting definitions leads to unintuitive guesswork and tedious verifications.

our approach

We believe that the best way to reckon with the technical challenges of hyperalgebra, while preserving its beneficial attributes, is to exploit a powerful relationship between powersets and complete lattices. Specifically, if $\mathcal{L}$ is a *suplattice*—a partially ordered set with all suprema—and $f : A \to \mathcal{L}$ a map of sets, then there is a unique extension of f to a map $f : \mathcal{P}(A) \to \mathcal{L}$ such that $f(\bigcup S) = \bigvee f(S)$. Similarly, we can uniquely lift lattice-valued operations to operations on $\mathcal{P}(A)$ that preserve joins in each variable. If we wish to recover the original operations, we can simply restrict to the singletons.

The aim of our project is to make use of this canonical lifting in order to study hyperstructures by studying their lattice counterparts. This incurs a cost: we take on the responsibility of managing an additional operation, the infinitary join $\bigvee$. However, by doing so we can retain the luxury of working with single-valued operations.

One may object to this approach on the grounds that it leads to an over-complicated picture of classical algebra. However, we have found the opposite to be the case. In particular, the theory of ideals of a commutative ring takes on a particularly elegant character when described in the language of lattices, a fact observed as early as the 1960s by Kirby.[8] And many of the well-known constructions from a first course in commutative algebra—the relationship between set-theoretic and ideal-theoretic product, the construction of colon ideals, and the prime ideal spectrum, to name three—can be seen to be general, universal constructions in this setting.

This project has already been carried out, albeit in a somewhat folkloric way, for the theory of semigroups and monoids. The corresponding lattice-enriched objects are called *quantales*, and we will need to establish their basic properties before turning our attention to our destination, lattice-enriched semirings.

QUANTALES AND POINTLESS TOPOLOGY

A *quantale* is, simply, a semigroup in the category of suplattices. They were first introduced by Mulvey in 1986 [11], where they were used to study quantum logic and non-commutative C*-algebras. The field of "pointless topology"—the work of Johnstone in particular—had used frames/locales (a certain kind of complete lattice) to give new insights into general topology, topos theory, spectra, and other things. Quantales, or "quantum locales", were supposed to capture non-commutative phenomena that go beyond the world described by frames.

But even commutative quantales already have an interesting theory. For example, the ideals in a commutative ring form a quantale under ideal-theoretic multiplication, and many old constructions in ring theory are most natural when seen from this perspective. An excellent treatment is given by Rosenthal in [14], from which we will borrow several technical results.

A *frame* is a quantale such that the multiplication equals the binary greatest lower bound. For our work, the most important thing about the category of quantales is that it reflects onto the category of frames. That is, for any quantale $\mathfrak{Q}$ there is a frame, which we suggestively call $\operatorname{Spec} \mathfrak{Q}$, and a morphism $\mathfrak{Q} \to \operatorname{Spec} \mathfrak{Q}$ initial among all morphisms of quantales $\mathfrak{Q} \to \mathcal{F}$, where $\mathcal{F}$ is a frame. In the general case, this construction is somewhat abstract, though it may be easy to write down in particular cases. For example, the universal frame associated to the quantale of ideals of a commutative ring is exactly the lattice of radical ideals.

Every frame $\mathcal{F}$ has an associated $T_0$ sober topological space $\operatorname{pt}(\mathcal{F})$, and a quotient map from $\mathcal{F}$ to the open sets of $\operatorname{pt}(\mathcal{F})$. As we can see, this map is sometimes an isomorphism, in which case we will say that $\mathcal{F}$ is *spatial*. Since $\operatorname{pt}(\mathcal{F})$ can actually be empty for nontrivial $\mathcal{F}$, it behooves us to regard $\operatorname{Spec} \mathfrak{Q}$ as a lattice, not a space. This is the perspective of *pointless topology*. Since the category of $T_0$ sober topological spaces is dual to the category of spatial frames and unital morphisms, this is a safe thing to do, provided that we can remember which way the arrows go.

One note: As written, Spec does not always produce unital morphisms, which means that we may wind up with too many morphisms. There is more than one way to fix this problem, but the easiest is just to restrict to the category of unital quantales and unital morphisms. This is an important point to make, since we will construct spectra for (hyper)semigroups and (hyper)semirings, but these constructions will only be functorial in a meaningful sense when we restrict to (hyper)monoids and unital (hyper)semirings.

Our basic method for writing down the spectrum of a semigroup $A$ is to first embed $A$ in its powerset $A \to \mathcal{P}(A)$, considered as a quantale, and then give some instructions for producing a quotient quantale $\mathcal{P}(A) \to \mathcal{C}(A)$. The second step is where we use the theory of quantales, while the first requires some digression about the power monad.

## motivation: the power monad

When Emmy Noether gave the first general definition of a commutative ring in 1921 [12], she was computing the algebras for a monad. Specifically, if $R$ denotes either $\mathbb{Z}$ or a field, then there is an endofunctor $T$ on **Set** sending $S$ to the set of polynomials $R[X_s]_{s \in S}$. There is a natural transformation $\eta : \mathrm{id} \to T$ sending $s$ to $X_s$, and another natural transformation $\mu : T^2 \to T$ given by expansion of polynomials. We have the following commutative diagrams:

$$\begin{array}{ccc}
T^3 \xrightarrow{T\mu} T^2 & T & T \xrightarrow{\eta T} T^2 \\
\mu T \downarrow \quad \downarrow \mu & T\eta \downarrow \searrow \mathrm{id} & \mathrm{id} \searrow \downarrow \mu \\
T^2 \xrightarrow{\mu} T & T^2 \xrightarrow{\mu} T & T
\end{array}$$

In modern terms, these diagrams say exactly that the functor $T$ is a *monad* over the category of sets, but establishing that they commute involves only basic properties of sets, functions, and polynomials, and is therefore entirely within the domain of 19[th] century mathematics.

Understanding the character of this monad, however, requires the modern definition of adjoint functor. The easiest way to create a monad from scratch is to compose a pair of adjoint functors, and all monads arise this way. However, a monad usually arises from many different adjunctions:

For any monad T on a category **C**, there is a category whose objects are categories **D** equipped with an adjunction **C** → **D** giving rise to T.

This category has both an initial object and a final object, which we can think of as the minimal adjunction and the maximal adjunction. The minimal adjunction is usually called the *Kleisli category* of T, and we can visualize it in two ways: as **C** itself, where a morphism A → B is a **C**-morphism A → TB, or as the category of the "free objects" TA, where morphisms TA → TB arise from lifts of A → TB. The maximal adjunction is the category of *Eilenberg-Moore algebras* for T, and it has many nice properties, characterized by Beck's monadicity theorem. For example, if **C** = **Set**, then the category of algebras admits kernel pairs and coequalizers.

The state of commutative algebra in the late 19$^{\text{th}}$ century was that there was tremendous interest in the Kleisli category of $T(S) = R[X_s]_{s \in S}$, the category of polynomial rings and polynomials maps. The great achievement of Noether's definition was to take us into the larger Eilenberg-Moore category of T-algebras, where the proofs are simpler and the examples richer.

The crucial monad in our work is the *power monad* sending a set S to its powerset $\mathcal{P}(S)$, where the unit $\eta : S \to \mathcal{P}(S)$ is given by $\eta(s) = \{s\}$, and the composition $\mu : \mathcal{P}^2(S) \to \mathcal{P}(S)$ is given by $\mu(\{S_\alpha\}) = \bigcup_\alpha S_\alpha$.

Hyperalgebra is stuck in the 19$^{\text{th}}$ century with respect to the power monad. That is, it takes place entirely within the Kleisli category of sets and hyperfunctions. Perhaps we will get simpler proofs, and a better theory, if we explore algebra instead in the Eilenberg-Moore category of $\mathcal{P}$: The category of suplattices. The philosophy of our work here is to apply this view fully to semigroups and semirings.

(Note: The general theory of lifting algebraic theories over a monad is beyond the scope of this work, but it seems essential to use 2-categories to explain the appearance of laxness.)

OUR RESULTS

After spending some time developing basic facts about lattices and quantales, we establish the following dictionary between classical algebraic structures, hyperalgebraic structures, and lattice-enriched ("quantic") algebraic structures:

| Classical Algebra | Hyperalgebra | Quantic Algebra |
|---|---|---|
| Sets | Sets | Suplattices |
| Semigroups | Hypersemigroups | Quantales |
| Monoids | Hypermonoids | Unital Quantales |
| Semirings | Hypersemirings | Quantic Semirings |

Of these, only quantic semirings are completely new, though we give a more general definition of hypersemiring than what currently appears in the literature, as we only assume subdistributivity, not strong distributivity.

Our main result is that, given a covered hypersemigroup A—that is, a hypersemigroup A equipped with a certain kind of closure operator—there is a universal morphism $D : A \to \operatorname{Spec} A$ (in the category of covered semigroups) to a frame. If the closure on A is finitary (generated by finite relations), then Spec A is the lattice of open sets of a $T_0$ sober topological space pt(Spec A), whose points are the prime ideals of A, and $\{D(a) \mid a \in A\}$ forms a basis of quasi-compact open sets. When restricted to covered hypermonoids, this gives a contravariant functor to the category of topological spaces.

Furthermore, we describe how to give a hypersemiring A a canonical structure of a finitary covered hypersemigroup under multiplication, giving a definition of the spectrum of a hypersemiring. This definition coincides with the usual definition of the (completely) prime spectrum of a (possibly non-commutative) ring, but it also coincides with the definition of Jun (and, earlier, Procesi-Ciampi and Rota [13]) of the spectrum of a Krasner hyperring.

Finally, we develop some connections of this theory to the ordered blueprints of Lorscheid.[10] An ordered blueprint satisfying a certain associative axiom is shown to give rise to a quantic semiring, which incidentally helps resolve an issue in the theory that the spectrum of a Krasner hyperring qua ordered blueprint is given too many points.

*Overview*

*Chapter 1:*

Chapter 1 is devoted to establishing the basic theory of closure operators on suplattices. In particular, we discuss the category of sets equipped with a closure operator on their powerset. We call these *closure domains*, instead of the more standard term *closure systems*, which we feel is more appropriately given to the collection of flats of a closure domain.

Whatever name they go by, closure domains are a very fundamental and general concept. For example, the category of topological spaces can be identified with the full subcategory of closure domains where the closure respects empty and binary unions. By placing other restrictions, we can identify other categories, such as matroids, convex geometries, and various kinds of lattices.

Suplattices and closures form the engine that makes this theory go. The closures on a suplattice are in bijection with its quotients, so they play a similar role as kernels do in the theory of abelian groups. In particular, the lattice of flats of a closure domain $A$ is naturally a quotient of the powerset lattice $(\mathcal{P}(A), \bigcup)$, which will be important for establishing universal properties later.

*Chapter 2:*

In Chapter 2, we look at semigroups and monoids in the category of closure domains, which we call *covered semigroups* and *covered monoids*, respectively. For a covered semigroup, the lattice of flats is no longer just a suplattice, but a *quantale*, and we spend some time discussing properties of such.

(Note: Many of the ideas in this chapter appear in a similar form in [15], though we take a slightly different focus—in particular, we take care to establish certain universal properties to more precisely apply the theory of quantales to our study of these objects.)

Since covered semigroups are in particular ordered semigroups, which are in particular categories, it is important to establish the essentially 2-categorical notions of *lax morphisms* and *colax morphisms* between them. (We do not directly use any higher category theory, but the connection is important and worth mentioning.) We also define $\mathbb{F}_1$ as a covered monoid, and show that the various kinds of morphisms to it are characterized by many familiar concepts, such as filters and ideals. In particular, prime ideals of a covered semigroup $A$ are in bijection with strong morphisms $A \to \mathbb{F}_1$. (In fact, these are the points of $\operatorname{Spec} A$, which we define later.)

We also look at the more general situation of a closure operator on a quantale that is compatible with the multiplication. These are known as *nuclei* in the literature, and they inherit many of the nice properties of closures on suplattices. In particular, they characterize quotient quantales. As such, the construction of certain nuclei is fundamental in showing the existence of universal quotient quantales, and we quote a few important results from Rosenthal on this matter.

In particular, every quantale has a universal quotient that is a *frame*—that is, a quantale where the multiplication coincides with the binary greatest lower bound. We use this to show that the spectrum of a commutative ring, as a covered monoid, has a universal morphism to a frame, namely the lattice of opens of its spectrum. In fact, we show that every covered semigroup A has a such a universal frame—by abuse of notation, we call this Spec A, identifying any sober $T_0$ space with its lattice of opens—which is necessarily the frame of opens of a topological space. One thing this gives us is a contravariant functor from the category of covered monoids to the category of spectral spaces.

*Interlude:*

Having spent Chapters 1 and 2 mostly giving new breath to old results, we turn, in the Interlude, towards something a bit new (though probably known, secretly, by the experts): the relationship between quantales and hypersemigroups.

In fact, we describe a general relationship between hyperoperations on a set A and lattice-enriched operations on the powerset lattice $(\mathcal{P}(A), \bigcup)$. We show that a hypersemigroup structure on a set A is the same as a quantale structure on $\mathcal{P}(A)$, and likewise a hypermonoid structure on A is the same as a unital quantale structure on $\mathcal{P}(A)$.

We use this to extend our results on spectra to covered hypersemigroups. We also show that if $\mathcal{X}$ is a topological space with a basis of quasi-compact opens, then the set A of all quasi-compact opens forms a covered hypersemigroup, with Spec A $= \mathcal{O}(\mathcal{X})$. In particular, if $\mathcal{X}$ is also sober and $T_0$, then Spec A $= \mathcal{X}$.

While associative binary hyperoperations lift nicely to $\mathcal{P}(A)$, we need a more delicate touch when dealing with more complicated hyperstructures, such as hyperrings.

*Chapter 3:*

Chapter 3 explores *quantic semirings*, which are designed to be the suplattice counterpart to hypersemirings. In fact, the definition that we obtain for hypersemirings is compatible with the definitions already given in the literature, and even somewhat more general, since we do not assume distributivity.

The key property of quantic semirings is *subdistributivity*:

$$S \cdot (U + T) \leqslant S \cdot U + S \cdot T$$

Unlike distributivity, subdistributivity lifts to the powerset.

Our key result is that quantales form a full reflective subcategory of quantic semirings. That is, there is a universal morphism from a quantic semiring to a quantale. For example, if A is a semiring, then the universal quantale associated to $\mathcal{P}(A)$ is the lattice of multiplicative sub-monoids of A. Since we can compute the spectrum of this quantale (and get the right answer!), this gives us a way to define the spectrum of a quantic semiring as a universal frame.

In general, the spectrum of a quantic semiring does not come from a topological space. However, we show that if A is a hypersemiring, then the spectrum of a quantic semiring $\mathcal{P}(A)$ coincides with the spectrum of A given a natural covered semigroup structure, and is therefore spatial. We show that this precisely recovers the spectrum of a hyperring constructed by Jun and others.

*Chapter 4:*

Chapter 4 explores examples and applications. One major application is to Lorscheid's theory of blueprints. We show how to interpret blue spectra using our language of covered monoids, and construct an additional nucleus, the nucleus of downward-closed sets, that gives additional insight into hyperrings as ordered blueprints. We show that, for any associative blueprint, the suplattice of downward-closed sets has natural operations that make it into a quantic semiring, whose spectrum is the space of ordered prime ideals. This suggests a new definition of spectra for ordered blueprints.

xiii

# CLOSURE SYSTEMS AND SUPLATTICES

## 1.1 CLOSURE SYSTEMS

**Definition 1.1.** A **closure system**[1] on a poset $\mathcal{P}$ is a subset $\mathcal{C} \subset \mathcal{P}$ such that, for every $S \in \mathcal{P}$, there is a minimum $T \in \mathcal{C}$ with $S \leqslant T$.

**Definition 1.2.** A **closure operator** on a poset $\mathcal{P}$ is a self-map $j : \mathcal{P} \to \mathcal{P}$ such that:

- $S \leqslant j(S)$.
- $S \leqslant T \implies j(S) \leqslant j(T)$.
- $j(j(S)) = j(S)$.

In the presence of the first two conditions, the third condition may be replaced by the weaker $j(j(S)) \leqslant j(S)$.

**Proposition 1.3.** *The following three structures on a poset $\mathcal{P}$ are equivalent:*

- *A closure system on $\mathcal{P}$.*
- *A closure operator on $\mathcal{P}$.*
- *A self-map $j : \mathcal{P} \to \mathcal{P}$ such that $S \leqslant j(T) \iff j(S) \leqslant j(T)$.*

*Proof.* The first and third conditions are easily seen to be equivalent, with $\mathcal{C} = j(\mathcal{P})$. If $j$ is a closure operator, and $S \leqslant j(T)$, then $j(S) \leqslant j(j(T)) = j(T)$. If $j(S) \leqslant j(T)$, then $S \leqslant j(S) \leqslant j(T)$. On the other hand, if $j$ satisfies the third condition, we check that it is a closure operator.

- $j(S) \leqslant j(S)$, so $S \leqslant j(S)$.
- If $S \leqslant T$, then $S \leqslant T \leqslant j(T)$, so $j(S) \leqslant j(T)$.
- $j(S) \leqslant j(S)$, so $j(j(S)) \leqslant j(S)$.

$\square$

---

[1] The term **Moore family** also appears in the literature.

*Remark* 1.4. The categorically-minded reader may notice that the third condition says exactly that $j : \mathcal{P} \to \mathcal{C}$ is left adjoint to the inclusion $\mathcal{C} \subset \mathcal{P}$. In fact, a closure operator on $\mathcal{P}$ is the same thing as a monad on $\mathcal{P}$, a fact which we do not use, but is nevertheless of great conceptual importance.

**Definition 1.5.** A **closure** on a poset $\mathcal{P}$ will mean the equivalent data of a closure system or a closure operator on $\mathcal{P}$. The elements of the associated closure system will be called **flats**, and, given a closure operator $j$, the collection of flats will be denoted $\mathcal{P}_j$.

*Remark* 1.6. We do not use the more common term "closed sets" for the elements of $\mathcal{C}$, as this might cause confusion later. In some cases, we will see that the flats of a closure will correspond to closed sets in a topology, but in cases of greater importance to us, they may correspond to open sets!

## 1.2 SUPLATTICES

**Definition 1.7.** A **suplattice** is a partially ordered set with all suprema. If $\mathcal{L}$ is a suplattice, and $\mathcal{S} \subset \mathcal{L}$ is a subset, we denote the supremum of $\mathcal{S}$ by $\bigvee \mathcal{S}$. Sometimes, if a subset of $\mathcal{L}$ is given as an indexed set $\{S_\alpha\}$, we will write the supremum as $\bigvee_\alpha S_\alpha$. A **morphism of suplattices** is a function $f : \mathcal{L} \to \mathcal{L}'$ that preserves the join: for any $\mathcal{S} \subset \mathcal{L}$, $f(\bigvee \mathcal{S}) = \bigvee f(\mathcal{S})$.

*Remark* 1.8. The partial order underlying a suplattice can be recovered from the join $\bigvee$, so it is often a good idea to present a suplattice as a set equipped with an infinitary operation $\bigvee : \mathcal{P}(\mathcal{L}) \to \mathcal{L}$ satisfying $\bigvee\{S\} = S$ for all $S \in \mathcal{L}$, and $\bigvee_\alpha \bigvee \mathcal{S}_\alpha = \bigvee \bigcup_\alpha \mathcal{S}_\alpha$ for all collections $\mathcal{S}_\alpha \subset \mathcal{L}$. (In this form, it is not hard to see that suplattices are exactly the algebras over the power monad, a perspective we will return to in the interlude.)

*Example* 1.9. If $S$ is a set, then the powerset $(\mathcal{P}(S), \subset)$, ordered by inclusion, is a suplattice, with join given by $\bigvee \mathcal{S} = \bigcup \mathcal{S}$.

**Proposition 1.10.** *On a partially ordered set, the existence of arbitrary least upper bounds is equivalent to the existence of arbitrary greatest lower bounds. In particular, a suplattice has all infima.*

*Proof.* If $\mathcal{L}$ is a suplattice, and $\mathcal{S} \subset \mathcal{L}$ is a subset, then let $\mathcal{T}$ be the set of all lower bounds for $\mathcal{S}$. It is immediate that $\bigvee \mathcal{T}$ is a greatest lower bound for $\mathcal{S}$. The reverse direction follows by the same argument applied to the dual poset. □

In particular, a suplattice has a binary meet, which we denote by $\wedge$. Note that the above argument does not require $\mathcal{T}$ to be nonempty, though this is in fact the case: every suplattice has a bottom element, the supremum of the empty set.

**Proposition 1.11.** *The set of closure operators on a suplattice $S$, ordered pointwise, is a suplattice.*

*Proof.* By Proposition 1.10, it is enough to show that it has infima. Given closure operators $j_\alpha$, in fact $j(a) = \bigwedge j_\alpha(a)$ is a closure operator. □

**Proposition 1.12.** *If $j$ is a closure operator on a suplattice $\mathcal{L}$, then $\mathcal{L}_j$ is also a suplattice.*

*Proof.* This follows from Proposition 1.10, and the fact that any infimum of flats is a flat. We could also show the existence of suprema more directly: If $\{C_\alpha\}$ is a collection of flats, then a flat $D$ contains each $C_\alpha$ if and only if $\bigvee_\alpha C_\alpha \subset D$. It follows that $j(\bigvee_\alpha C_\alpha)$ is the least upper bound for $\{C_\alpha\}$ in $\mathcal{L}_j$. □

Note that $j : \mathcal{L} \to \mathcal{L}_j$ is a morphism of suplattices. We could also characterize the suplattice structure on $\mathcal{L}_j$ as the unique one making $j$ into a morphism of suplattices. In fact, $\mathcal{L}_j$ is a quotient of $\mathcal{L}$, in the following sense:

**Proposition 1.13.** *Suppose that $\mathcal{L}, \mathcal{L}'$ are suplattices, $j$ a closure operator on $\mathcal{L}$. If $f : \mathcal{L} \to \mathcal{L}'$ is a morphism of suplattices with $f(S) = f(j(S))$ for all $S \in \mathcal{L}$, then there is a unique morphism $f_j : \mathcal{L}_j \to \mathcal{L}'$ such that the following diagram commutes:*

$$\begin{array}{ccc} \mathcal{L} & \xrightarrow{j} & \mathcal{L}_j \\ & \searrow{\scriptstyle f} & \downarrow{\scriptstyle f_j} \\ & & \mathcal{L}' \end{array}$$

*Proof.* It is clear that $f_j$ exists and is unique as a map of sets—we must have $f_j(C) = f(C)$—so all we need to check is that it is really a morphism of suplattices. If $\bigvee_j$ denotes the join in $\mathcal{L}_j$, we have:

$$f_j(\bigvee_j S) = f(\bigvee_j S) = f(j(\bigvee S)) = f(\bigvee S) = \bigvee f(S) = \bigvee f_j(S)$$

□

Another way to see that closures characterize quotients is to compare them to congruence relations:

**Definition 1.14.** A **congruence relation** on a suplattice $\mathcal{L}$ is an equivalence relation $\sim$ such that, if $s \sim t$ for all $s \in S, t \in T$, then $\bigvee S \sim \bigvee T$.

**Proposition 1.15.** *The closure operators and congruence relations on a suplattice are in bijection.*

*Proof.* Given a congruence relation, we can define a closure operator $j(S) = \bigvee_{S \sim T} T$. Conversely, a closure operator gives a congruence: $S \sim T$ if $j(S) = j(T)$. □

## 1.3 CLOSURE DOMAINS

**Definition 1.16.** A **closure domain** is a set a equipped with a closure on $(\mathcal{P}(A), \subset)$. We will also refer to a **closure operator** or **closure system** on $A$, when this will not cause confusion. The lattice of flats of a closure domain $A$ will be denoted $\mathcal{C}(A)$.

*Remark* 1.17. We will generally use the notation $S \mapsto [S]$ to describe closure operators that live on powersets. Having a second notation will be useful when both a poset and its powerset possess a closure operator.

**Definition 1.18.** A **closure relation** on $A$ is a relation $\vdash$ between $\mathcal{P}(A)$ and $A$, such that:

- $a \in S \implies S \vdash a$.

- If $T \vdash s$ for all $s \in S$, and $S \vdash a$, then $T \vdash a$.

*Remark* 1.19. It may be helpful to pronounce $S \vdash a$ as ""S covers $a$", and to call $\vdash$ a "covering relation". However, we will prefer to reserve the terminology of covers for the situation of Chapter 2, where the closure is compatible with some semigroup operation.

**Proposition 1.20.** *Closure relations on $A$ are in bijection with closure systems and closure operators on $A$.*

*Proof.* Given a closure relation $\vdash$, define $[S] = \{a \mid S \vdash a\}$. Conversely, given a closure operator $[-]$, define $S \vdash a$ if $a \in [S]$. We omit the verification that these give a bijection between closure relations and closure operators. □

**Proposition 1.21.** *Let $A$ and $B$ be closure domains, $f : A \to B$ a function on the underlying sets. The following three conditions are equivalent:*

- *For all $C \in \mathcal{C}(B)$, $f^{-1}(C) \in \mathcal{C}(A)$.*

- *For all $S \subset A$, $f([S]) \subset [f(S)]$.*

- *For all $S \subset A$, $a \in A$, $S \vdash a \implies f(S) \vdash f(a)$.*

*Proof.* The second and third conditions are easily seen to be equivalent.

If the first condition is satisfied, then, since $[f(S)]$ is a flat, $f^{-1}([f(S)])$ is also a flat. We have $S \subset f^{-1}(f(S)) \subset f^{-1}([f(S)])$, so $[S] \subset f^{-1}([f(S)])$, therefore $f([S]) \subset [f(S)]$.

Suppose the second condition is satisfied. Let C be a flat, and set $S = f^{-1}(C)$. Since $f([S]) \subset [f(S)] = [C] = C$, we have $[S] \subset f^{-1}(f([S])) \subset f^{-1}(C) = S$. Since $S \subset [S]$, we have $[S] = S$ and S is a flat. □

**Definition 1.22.** A **morphism of closure domains** is a function $f : A \to B$ satisfying any and all of the above three conditions. A **covering morphism** is a morphism such that $[f(A)] = B$.

*Remark* 1.23. We include the notion of covering morphism because it may have algebraic significance in the absence of a multiplicative unity—for example, a (not necessarily unital) morphism $f : A \to B$ of commutative rings with unity is unital if and only if $f(A)$ generates B as an ideal. For the most part, we will focus on the unital case for simplicity, and briefly discuss, in an aside, the use of covering morphisms in the non-unital case. The main thing to note, for now, is that identity maps are covering morphisms, and covering morphisms and are closed under composition, and therefore covering morphisms define a subcategory.

**Definition 1.24.** A closure or closure domain is called **finitary** or **algebraic** if, whenever $S \vdash a$, there is a finite subset $S' \subset S$ such that $S' \vdash a$. Equivalently:

$$[S] = \bigcup_{\substack{S' \subset S \\ S' \text{ finite}}} [S']$$

## 1.4 EXAMPLES OF CLOSURE DOMAINS

*Algebraic examples*

In general, set-theoretic models for algebraic theories give rise to finitary closure domains.[2] Here are some examples:

- A is a group, and [S] is the subgroup generated by S.

---

2 Here we use the word "model" in a technical sense, for instance a group is a set-theoretic model of the first-order theory of groups.

- A is a left R-module for a ring R, and [S] is the submodule of A generated by S.

- A is a left G-set for a group G, and [S] is the the union of the orbits meeting S.

These examples all arise from monads: If $\mathcal{T}$ is a monad on **Set**, and A is a $\mathcal{T}$-algebra with structure map $\mathcal{T}A \to A$, then we can define a closure on A by taking [S] to be the image of $\mathcal{T}S \to \mathcal{T}A \to A$.

But there are also interesting examples that do not appear to arise from monads on **Set**:

- A is a group, and [S] is the normal closure of the subgroup generated by S.

- $\mathcal{P}(X)$ is a powerset, and, for some collection $\mathcal{S} \subset \mathcal{P}(X)$, [$\mathcal{S}$] is the coarsest topology containing $\mathcal{S}$.

*Galois connections*

If $\mathcal{P}, \mathcal{Q}$ are posets, and $f : \mathcal{P} \to \mathcal{Q}, g : \mathcal{Q} \to \mathcal{P}$ is a Galois connection (either covariant or contravariant), then $g \circ f$ is a closure operator on $\mathcal{P}$. In fact, all closures arise this way, as we will see later.

*The underlying preorder*

A **preordered set** is a set P equipped with a reflexive and transitive relation $\leqslant$. A morphism $f : P \to P'$ of preordered sets is a function such that $a \leqslant b \implies f(a) \leqslant f(b)$.

Preordered sets embed as a full subcategory of closure domains: the closure relation is given by $S \vdash a$ if $a \leqslant s$ for some $s \in S$. The flats of this closure are exactly the sets that are downward-closed: C is a flat exaclty if $x \leqslant y$ and $y \in C$ implies $x \in C$.

A closure domain arises from a preordered set if and only if, whenever $S \vdash a$, there is some $s \in S$ such that $\{s\} \vdash a$. And given a closure domain A, we have a natural preorder on A, given by $a \leqslant b$ if $a \in [\{b\}]$. In fact, this is a coreflection:

**Proposition 1.25.** *Let* P *be a preordered set,* A *a closure domain. A function* $f : P \to A$ *is a morphism of preorders if and only if it is a morphism of closure domains.*

*Proof.* Suppose that f is a morphism of preorders, and $S \vdash a$. Then there is some $s \in S$ with $a \leqslant s$, so $f(a) \leqslant f(s)$. Then $[\{f(s)\}] \vdash f(a)$, so $f(S) \vdash f(a)$.

Conversely, if f is a morphism of closure systems, and $a \leqslant b$, then $\{b\} \vdash a$, so $\{f(b)\} \vdash f(a)$, and $f(a) \leqslant f(b)$. □

**Definition 1.26.** A closure or closure domain is called **proper** if its underlying preorder is a partial order.

*Topological Spaces I*

Let X be a topological space. Then $S \mapsto \overline{S}$ is a closure on X, with flats the closed sets. Of course, the topology on X can be recovered directly from the closure system.

In fact, by Proposition 1.21, this gives a full and faithful embedding of the category of topological spaces as a full subcategory of the category of closure domains. By an old insight of Kuratowski [9], we can characterize its essential image as those closure domains satisfying the following two properties:

- $[\emptyset] = \emptyset$.

- $[S \cup T] = [S] \cup [T]$.

These two axioms, together with those for a closure operator, are frequently named in the literature as the **Kuratowski closure axioms**.

*Topological spaces II*

Let $\mathfrak{X}$ be a topological space, and let B be a basis. Then there is a natural closure on B, given by $S \vdash a$ if $a \subset \bigcup S$, where $a \in B$ is an element of the basis and $S \subset B$ is a collection of basis elements. The flats with respect to this closure are the open sets in $\mathfrak{X}$. In particular, the collection of all open sets, $\mathcal{O}(\mathfrak{X})$, is a closure domain.[3]

Note that if $f : \mathfrak{X} \to \mathfrak{Y}$ is a continuous map of topological spaces, the natural morphism of closure domains goes in the opposite direction: $f^{-1} : \mathcal{O}(\mathfrak{Y}) \to \mathcal{O}(\mathfrak{X})$.

This example reflects what we intend to do with closures much more than the previous one. We are primarily interested in algebraic examples with geometric flavor, and so our morphisms will generally point in the opposite direction of the natural geometric interpretation. This will be revisted in the next chapter.

---
[3] In fact it is a suplattice, which will be discussed later in the section, and furthermore a frame, which will be discussed in the next chapter.

*Join-semilattices*

A **join-semilattice** is a set L equipped with an associative, commutative, idempotent binary operation $\vee$. Equivalently, L is a poset such that any pair of elements has a least upper bound. A morphism $f : L \to L'$ of join-semilattices is a function such that $f(a \vee b) = f(a) \vee f(b)$ for all $a, b \in L$.

A join-semilattice L comes equipped with a canonical finitary closure: $S \vdash a$ if $a \leqslant s_1 \vee \cdots \vee s_n$ for some $s_i \in S$. The flats of this closure are usually referred to in lattice theory as *ideals*: they are the sets C that are downward-closed, and satisfy $x, y \in C \implies x \vee y \in C$.

If L and L' are join-semilattices, then a function $f : L \to L'$ is a morphism of join-semilattices if and only if it is a morphism of closure domains. So join-semilattices are a full subcategory of closure domains.

*Directed-complete partial orders*

A **directed-complete partial order** or **dcpo** is a poset P such that every directed subset has a least upper bound. Morphisms of dcpos are functions preserving the directed join, which are often called *Scott-continuous* functions.

Once again, these form a full subcategory of closure domains. The closure on a dcpo is given by: $S \vdash a$ if there is a directed set $D \subset S$ with $a \leqslant \bigvee D$. Its flats are the Scott-closed sets: the downward-closed sets closed under directed joins.

*Matroids*

A **matroid** is a closure domain satisfying the **exchange axiom**: $a \in [S \cup \{b\}] \setminus [S] \implies b \in [S \cup \{a\}]$. A wide class of matroids arise from vector spaces:

Let k be a field, V a vector space over k, and take any subset $A \subset V$. Then A is naturally a closure domain, with closure given by $[S] = A \cap \langle S \rangle$, where $\langle S \rangle$ denotes the vector space generated by S. By standard linear algebra, A is furthermore a matroid.

The morphisms of matroids that are morphisms of closure domains are called **strong morphisms** by some authors.

*Convex geometries and antimatroids*

A **convex geometry** is a closure domain satisfying $[\emptyset] = \emptyset$, and the **anti-exchange axiom**: if $a \neq b$, and $a, b \notin [S]$, then $a \in [S \cup \{b\}] \implies b \notin [S \cup \{a\}]$. The flats of a convex geometry are called **convex sets**.

If $\mathbb{K}$ is an ordered field, and $A \subset \mathbb{K}^n$, then there is a closure on $A$ taking $S$ to the intersection of $A$ with the convex hull of $S$. A basic result of convex geometry is that this operator is finitary: in fact, Carathéodory's Theorem says that if $a \in [S]$, then $a \in [S']$ for some $S'$ with $|S'| \leqslant n+1$—the convex hull of $S$ is the union of the (possibly degenerate) n-simplices with vertices in $S$. Since simplices are exactly the finitely generated flats, $S \mapsto [S]$ is finitary.

An **antimatroid** is a set $A$ and a collection $\mathcal{D}$ of subsets of $A$ whose complements form a convex geometry. Antimatroids themselves are not literally closure domains, but they can be thought of as coclosures: closures on the dual poset $(\mathcal{P}(A), \supset)$.

There does not seem to be a widely accepted definition of morphisms for convex geometries and antimatroids, but the morphisms as closure domains are certainly of interest. For example, a function $f : \mathbb{K}^m \to \mathbb{K}^n$ is a morphism of closure systems if and only if $f$ preserves *betweenness*: if $x$ lies on the line segment between $p$ and $q$, then $f(x)$ lies on the line segment between $f(p)$ and $f(q)$. Affine transformations have this property, but characterizing all such functions appears to be a difficult problem, even for $\mathbb{K} = \mathbb{R}$.

*Closures on a ring*

Let $A$ be a commutative ring with 1. $A$ comes equipped with many closures of interest. We enumerate, suggestively, some of the notable ones:

- $S \mapsto S$.

- $S \mapsto \lfloor S \rfloor$, the (additive) semigroup generated by $S$.

- $S \mapsto [S]$, the (additive) monoid generated by $S$.

- $S \mapsto \langle S \rangle$, the (additive) abelian group generated by $S$.

- $S \mapsto (S)$, the ideal generated by $S$.

- $S \mapsto |S|$, the radical ideal generated by $S$.

We have $S \subset [S] \subset \langle S \rangle \subset (S) \subset |S|$ for all $S$, so each of these closures factors through the previous one. Each subsequent extension loses information, but gains algebraic structure:

- $\lfloor S \rfloor + \lfloor S \rfloor \subset \lfloor S \rfloor$.
- $0 \in [S]$.
- $-1 \cdot \langle S \rangle \subset \langle S \rangle$.
- $R \cdot (S) \subset (S)$.
- $x^2 \in |S| \implies x \in |S|$.

Critically, in the final step we have a closure whose flats are the open sets of a topological space $\mathrm{Spec}\,A$. As we proceed, we will see that this is not at all a coincidence: we will describe the general procedure by which closures may be employed to turn algebra into topology.

## 1.5 SUPLATTICES AS CLOSURE DOMAINS

If $\mathcal{L}$ is a suplattice, we can consider it as a closure domain, with $[S] = \{a \mid a \leqslant \bigvee S\}$. The flats are the principal downward-closed sets.

This gives us an equivalence between the category of suplattices and the full subcategory of closure domains such that, for any $S \subset A$, there is a unique $a \in A$ such that $[S] = [\{a\}]$. The next few propositions will show that, in fact, this category is reflective.

**Proposition 1.27.** *If $A$ is a closure domain, then $\mathcal{C}(A)$ is a suplattice under inclusion, and the natural map $A \to \mathcal{C}(A)$ is a morphism of closure domains.*

*Proof.* The first part follows from Proposition 1.12. To check that the map $\iota \colon A \to \mathcal{C}(A)$ sending $a$ to $[\{a\}]$ is a morphism of closure domains, suppose that $S \vdash a$. Then $a \in [S]$, so $[\{a\}] \subset [S] = \bigvee_{s \in S}[\{s\}]$. So $\iota(a) \subset \bigvee \iota(S)$, therefore $\iota(S) \vdash \iota(a)$. □

**Proposition 1.28.** *Let $A$ be a set, $\mathcal{L}$ a suplattice. If $f \colon A \to \mathcal{L}$ is any set map, there is a unique morphism $\mathcal{P}(f) \colon \mathcal{P}(A) \to \mathcal{L}$ of suplattices making the following diagram commute:*

$$\begin{array}{ccc} A & \longrightarrow & \mathcal{P}(A) \\ & \searrow{\scriptstyle f} & \downarrow{\scriptstyle \mathcal{P}(f)} \\ & & \mathcal{L} \end{array}$$

*Proof.* If $\mathcal{P}(f)$ is such a morphism of suplattices, then we have:

$$(\mathcal{P}(f))(S) = (\mathcal{P}(f))(\bigcup S) = \bigvee_{s \in S}(\mathcal{P}(f))(\{s\})$$

$$= \bigvee_{s \in S} f(s) = \bigvee f(S)$$

It remains to check that $(\mathcal{P}(f))(S) = \bigvee f(S)$ defines a morphism of suplattices. (It obviously makes the diagram commute.)

$$(\mathcal{P}(f))(\bigcup_\alpha S_\alpha) = \bigvee \bigcup_\alpha S_\alpha$$

$$= \bigvee_\alpha \bigvee S_\alpha = \bigvee_\alpha (\mathcal{P}(f))(S_\alpha)$$

$\square$

**Proposition 1.29.** *Let $A$ be a closure domain, $\mathcal{L}$ a suplattice. If $f : A \to \mathcal{L}$ is a morphism of closure domains, then there is a unique morphism of suplattices $\mathcal{C}(A) \to \mathcal{L}$ making the following diagram commute:*

$$\begin{array}{ccc} A & \longrightarrow & \mathcal{C}(A) \\ & \searrow{\scriptstyle f} & \downarrow{\scriptstyle \exists!} \\ & & \mathcal{L} \end{array}$$

*Proof.* This is a consequence of Proposition 1.28 and Proposition 1.13: first we lift $f$ to $\mathcal{P}(f) : \mathcal{P}(A) \to \mathcal{L}$, then to a morphism $\mathcal{C}(A) = [\mathcal{P}(A)] \to \mathcal{S}$. The only thing we need to check is that $\mathcal{P}(f)$ respects the closure. But this is equivalent to the statement that $f$ is a morphism of closure domains. $\square$

# 2

# QUANTALES, $\mathbb{F}_1$, AND THE POINTLESS SPECTRUM

## 2.1 COVERED SEMIGROUPS

Let $(A, \cdot)$ be a semigroup. By convention, we extend the multiplication to the powerset $\mathcal{P}(A)$ by setting $S \cdot T = \{s \cdot t \mid s \in S, t \in T\}$.

**Definition 2.1.** A **coverage** on $(A, \cdot)$ is a closure $[-]$ that satisfies:

$$[S] \cdot [T] \subset [S \cdot T] \text{ for all } S, T \subset A$$

Equivalently:

$$a \cdot [T] \subset [a \cdot T] \text{ and } [S] \cdot a \subset [S \cdot a] \text{ for all } a \in A, S, T \subset A$$

A **covered semigroup** is a semigroup equipped with a coverage.

*Remark* 2.2. To see that these two conditions are equivalent, notice that the second is the same as the assertion that $U \cdot [V] \subset [U \cdot V]$ and $[U] \cdot V \subset [U \cdot V]$ for all $U, V \subset A$. Combining these gives us:

$$[S] \cdot [T] \subset [S \cdot [T]] \subset [[S \cdot T]] = [S \cdot]$$

*Remark* 2.3. In the case that $(A, [-])$ is a partially ordered set, this is identical to the usual notion of an ordered semigroup. In fact, the following proposition shows that this definition is compatible with the forgetful functor to preordered semigroups.

**Proposition 2.4.** *If $A$ is a covered semigroup, then $\cdot$ is monotone on $(A, \leqslant)$.*

*Proof.* If $a \leqslant b$ and $c \leqslant d$, then $a \in [b]$ and $c \in [d]$, so $a \cdot c \in [b] \cdot [d] \subset [b \cdot d]$, and therefore $a \cdot c \leqslant b \cdot d$. □

**Definition 2.5.** If $A$ is a covered semigroup, a **left ideal** of $A$ is a flat $I$ such that $A \cdot I \subset I$, and a **right ideal** is a flat $I$ such that $I \cdot A \subset I$. A **ideal** is a flat that is both a left ideal and a right ideal. A **prime ideal** is an ideal such $I \neq A$, and $x \cdot y \in I \implies x \in I$ or $y \in I$. A **prime filter** is the complement of a prime ideal.

*Remark* 2.6. In the case of non-commutative multiplication, it may be more familiar to use the name "completely prime ideal" instead of "prime ideal". We use the simpler term to emphasize a more unified approach, though we will see presently that these definitions are more well-behaved in the commutative case.

**Proposition 2.7.** *If $A$ is a commutative covered semigroup, there is a smallest coverage $(-)$ above $[-]$ such that $A \cdot (S) \subset (S)$, and the flats of this coverage are exactly the ideals of $A$.*

*Proof.* Since $A \cdot [S \cup A \cdot S] \subset [A \cdot S \cup A \cdot A \cdot S] = [A \cdot S] \subset [S \cup A \cdot S]$, it is enough to show that $(S) = [S \cup A \cdot S]$ is a coverage. First, we check that it is a closure operator:

- $S \subset S \cup A \cdot S \implies S \subset (S)$.

- $S \subset T \implies S \cup A \cdot S \subset T \cup A \cdot T \implies (S) \subset (T)$.

- $((S)) = [(S) \cup A \cdot (S)] = [(S)] = [[S \cup A \cdot S]] = [S \cup A \cdot S] = (S)$.

Next, we check multiplicativity:

$$(S) \cdot (T) = [S \cup A \cdot S] \cdot [T \cup A \cdot T]$$
$$\subset [(S \cup A \cdot S) \cdot (T \cup A \cdot T)]$$
$$= [S \cdot T \cup A \cdot S \cdot T] = (S \cdot T)$$

Finally, if $S$ is an ideal of $A$, we have $(S) = [S \cup A \cdot S] = [S] = S$, so $S$ is a flat of $(-)$, and conversely, if $S$ is a flat of $(-)$, $S = (S) = [S \cup A \cdot S] \implies A \cdot S \subset S$, and $S$ is an ideal of $A$. □

*Remark* 2.8. In the non-commutative case, the above proposition becomes much more complicated. It is still true that there is a smallest coverage $(-)$ such that every flat is a left ideal, but it is no longer true that the left ideals of $[-]$ are necessarily flats of $(-)$. This inconvenient fact is responsible for the divergence of the concepts of prime ideal and completely prime ideal in the theory of noncommutative rings, and it means that great care should be taken in applying our results in the noncommutative setting.

For more details on this construction, see Chapter 3 of [14], or our outline of the spectrum construction later in the chapter.

**Definition 2.9.** A **covered monoid** is a covered semigroup with a multiplicative unit.

*Remark* 2.10. If A is a covered monoid, the construction of Proposition 2.1 can be simplified: we can take $(S) = [A \cdot S]$.

**Definition 2.11.** A covered monoid is **affine** or **strictly two-sided** if 1 is a top element of its preorder.

## 2.2 SEMILATTICES

**Definition 2.12.** A is a **closure domain with meets** if the underlying preorder is a meet-semilattice. That is, A is proper, and every pair of elements has a greatest lower bound. Likewise, a **covered semigroup with meets** (resp. **covered monoid with meets**) is a covered semigroup (resp. covered monoid) such that the underlying closure domain is a closure domain with meets.

**Definition 2.13.** A **covered semilattice** is a proper covered semigroup such that $a \cdot b$ is a greatest lower bound for $\{a, b\}$ for all $a, b \in A$. Equivalently, it is a closure domain whose underlying preorder is a meet-semilattice, and $[S] \wedge [T] \subset [S \wedge T]$ for all $S, T \subset A$.

*Remark* 2.14. It is worthwhile to note that these two equivalent definitions mean that we can think of covered semilattices as either a full subcategory of covered semigroups, or a non-full subcategory of closure domains.

*Remark* 2.15. Coverages on semilattices are closely related to the structures of the same name studied by Johnstone. (see, e.g. Section II.2.11 of [6].) The differences are mainly philosophical: For example, Johnstone does not require $\{s\}$ to be a cover of $s$. A meaningful comparison is that his coverages are similar to Grothendieck pretopologies, while ours are similar to Grothendieck topologies.

**Proposition 2.16.** *Suppose that* $(A, \cdot, 1)$ *is a covered monoid. Then A is a covered semilattice if and only if it is proper, idempotent, and affine.*

*Proof.* In a covered semilattice with 1, it is easy to see that $\wedge$ is idempotent and 1 is the top element. Conversely, A is proper, idempotent, and affine, then we have, for all $a, b \in A$:

$$a \cdot b \leqslant a \cdot 1 = a$$
$$a \cdot b \leqslant 1 \cdot b = b$$

So $a \cdot b$ is a lower bound for $\{a, b\}$. On the other hand, if $x$ is any lower bound of $\{a, b\}$, we have:

$$x = x \cdot x \leqslant a \cdot b$$

It follows that $a \cdot b$ is a greatest lower bound of $\{a, b\}$. Since A is proper, we have $\cdot = \wedge$ and A is a covered semilattice. □

## 2.3 LAX AND COLAX MORPHISMS

**Definition 2.17.** If $A, B$ are covered semigroups, a **weak morphism** is just a morphism of closure domains $f : A \to B$. We say that a weak morphism $f$ is a **lax morphism** if it satisfies:

$$f(a) \cdot f(b) \leqslant f(a \cdot b)$$

Likewise, $f$ is a **colax morphism** if it satisfies:

$$f(a \cdot b) \leqslant f(a) \cdot f(b)$$

$f$ is a **strong morphism** if it is both lax and colax. $f$ is a **strict morphism** if it is a morphism of semigroups in the usual sense. Note that strict morphisms are strong, and, if B is proper, strong morphisms are strict.

*Remark* 2.18. There are many examples in the literature that use the reverse of our definitions for lax and colax. We have good reasons for choosing this terminology, however: First, we wish for an inequality $a \leqslant b$ to be treated as a morphism $a \to b$ rather than $b \to a$, to respect the convention that an inclusion of sets is a morphism. Second, the notions of lax and colax appear in the theory of monoidal categories, and we wish to be consistent with the terminology there.

**Proposition 2.19.** *If A is a closure domain and semigroup, then it is a covered semigroup if and only if $[-]$ is a lax morphism $\mathcal{P}(A) \to \mathcal{P}(A)$.*

*Proof.* This is a formal consequence of the definitions, since A is a covered semigroup if and only if $[S] \cdot [T] \subset [S \cdot T]$. □

**Definition 2.20.** If A and B are covered monoids, we call any of the above classes of morphism **unital** if $f(1) = 1$.

*Remark* 2.21. If A is a monoid, B is a covered monoid, and $f : A \to B$ is a morphism of closure domains, we can define *lax-unital* morphisms to be those with $1 \leqslant f(1)$, and *colax-unital* morphisms to be those with $f(1) \leqslant 1$. We can similarly define the terms *strong-unital* and *strict-unital*.

Here we will assume strictness for simplicity—though some results generalize to lax-unital morphisms—but there are interesting examples in the literature of various combinations of laxness and colaxness for the multiplication and unit. For example, [3] considers *lazy morphisms*, which are colax and lax-unital.

## 2.4 QUANTALES

**Definition 2.22.** A **quantale** is a triple $(\mathcal{Q}, \bigvee, \cdot)$, where $(\mathcal{Q}, \bigvee)$ is a suplattice and $(\mathcal{Q}, \cdot)$ is a semigroup, that satisfies, for all $U, V, U_\alpha, V_\alpha \in \mathcal{Q}$:

$$U \cdot \left(\bigvee_\alpha V_\alpha\right) = \bigvee_\alpha (U \cdot V_\alpha)$$

$$\left(\bigvee_\alpha U_\alpha\right) \cdot V = \bigvee_\alpha (U_\alpha \cdot V)$$

A **unital quantale** is a quantale with a unit. A (lax/colax/strong) morphism of quantales is a function $f : \mathcal{Q} \to \mathcal{Q}'$ that is both a morphism of suplattices and a (lax/colax/strong) morphism of ordered semigroups.

In most cases, it is probably best to think of quantales concretely as sets equipped with two operations, $\cdot : \mathcal{Q}^2 \to \mathcal{Q}$ and $\bigvee : \mathcal{P}(\mathcal{Q}) \to \mathcal{Q}$, or as partially ordered sets with suprema and a compatible multiplication.

But crucially, quantales (with lax/colax/strong morphisms) can be realized as a full subcategory of covered semigroups (with lax/colax/strong morphisms). Specifically, they are the covered semigroups whose underlying closure system can be identified with a suplattice. Put differently, they are the covered semigroups such that every flat is principal with a unique generator.

*Example* 2.23. If $A$ is a covered semigroup, then the multiplication $C \cdot_{[-]} D = [C \cdot D]$ makes $\mathcal{C}(A)$ into a quantale. In particular, if $(A, \cdot)$ is a semigroup, then $(\mathcal{P}(A), \bigcup, \cdot)$ is a quantale. (We will see later that $A \to \mathcal{C}(A)$ is a reflection.)

**Definition 2.24.** A **nucleus** on a quantale $\mathcal{Q}$ is a closure $j$ on the underlying suplattice, such that $j(S) \cdot j(T) \leqslant j(S \cdot T)$ for all $S, T \in \mathcal{Q}$.

By comparing definitions, we can see that a coverage on a semigroup $(A, \cdot)$ is precisely the same as a nucleus on the quantale $(\mathcal{P}(A), \bigcup, \cdot)$. We will tend to use the words "coverage" and "nucleus" somewhat interchangeably, the former when we wish to emphasize the underlying set $A$, and the latter when we wish to emphasize the lattice of subsets $(\mathcal{P}(A), \bigcup)$.

As with closure operators on suplattices, nuclei characterize quotients of quantales. In particular, they have the following universal property:

**Proposition 2.25.** *Suppose that $\mathcal{Q}, \mathcal{Q}'$ are quantales, $[-]$ a nucleus on $\mathcal{Q}$. If $f : \mathcal{Q} \to \mathcal{Q}'$ is a lax(/colax) morphism of quantales with $f(S) = f([S])$ for all $S \in \mathcal{Q}$, then there is a unique lax(/colax) morphism $[f] : [\mathcal{Q}] \to \mathcal{Q}'$ such that the following diagram commutes:*

$$\begin{array}{ccc} \mathcal{Q} & \longrightarrow & [\mathcal{Q}] \\ & \searrow_{f} & \downarrow_{[f]} \\ & & \mathcal{Q}' \end{array}$$

*Proof.* Since $[-]$ is a closure, we can use the universal property of closures on suplattices to construct $[f]$ as a morphism of suplattices: $[f](C) = f(C)$. It remains only to check that $[f]$ is lax(/colax). If $f$ is lax:

$$[f](S) \cdot [f](T) = f(S) \cdot f(T) \leqslant f(S \cdot T) = f([S \cdot T]) = [f]([S \cdot T]) = [f](S \cdot_{[-]} T)$$

If $f$ is colax, the same proof works with the inequality reversed. $\square$

**Proposition 2.26.** *Let $A$ be a semigroup, $\mathcal{Q}$ a quantale. If $f : A \to \mathcal{Q}$ is a lax(/colax) morphism, then there is a unique lax(/colax) morphism of quantales $\mathcal{P}(f) : \mathcal{P}(A) \to \mathcal{Q}$ making the following diagram commute:*

$$\begin{array}{ccc} A & \longrightarrow & \mathcal{P}(A) \\ & \searrow_{f} & \downarrow_{\mathcal{P}(f)} \\ & & \mathcal{Q} \end{array}$$

*Proof.* By the earlier universal property of the powerset, it is enough to show that the natural map $\mathcal{P}(f)$ is lax(/colax) when $f$ is. If $f$ is lax:

$$(\mathcal{P}(f))(S \cdot T) = \bigvee_{s \in S, t \in T} f(s \cdot t) \leqslant \bigvee_{s \in S, t \in T} f(s) \cdot f(t)$$

$$= (\bigvee_{s \in S} f(s)) \cdot (\bigvee_{t \in T} f(t)) = ((\mathcal{P}(f))(S)) \cdot ((\mathcal{P}(f))(T))$$

If $f$ is colax, we just reverse the inequality. $\square$

**Proposition 2.27.** *Let $A$ be a covered semigroup, $\mathcal{Q}$ a quantale. If $f : A \to \mathcal{Q}$ is a lax(/colax) morphism of covered semigroups, then there is a unique lax(/colax) morphism of quantales $\mathcal{C}(A) \to \mathcal{Q}$ making the following diagram commute:*

$$\begin{array}{ccc} A & \longrightarrow & \mathcal{C}(A) \\ & \searrow{\scriptstyle f} & \downarrow{\scriptstyle \exists!} \\ & & \mathcal{Q} \end{array}$$

*Proof.* As in Chapter 1, this is an immediate consequence of the previous two propositions. □

*Remark* 2.28. Once we have developed the theory of frames, this last proposition can be seen as a generalization of Stone duality. For example, it immediately gives us an equivalence of categories between bounded distributive lattices and coherent frames with coherent maps, the latter of which is easily seen to be dual to the category of spectral spaces.

## 2.5 FRAMES

Recall that a suplattice has infima, in particular it has a binary meet $\wedge$.

**Definition 2.29.** A **frame** is a suplattice $(\mathcal{F}, \bigvee)$ such that, for all $U, V_\alpha \in \mathcal{F}$, we have:

$$U \wedge \left( \bigvee_\alpha V_\alpha \right) = \bigvee_\alpha (U \wedge V_\alpha)$$

A morphism of frames is a function $f : \mathcal{F} \to \mathcal{F}'$ that is a morphism of suplattices, and satisfies $f(x \wedge y) = f(x) \wedge f(y)$.

The category of frames, as we have defined it, is easily seen to be a full subcategory of the category of quantales with strong morphisms. We have found it useful to think of frames as the covered semigroups which are simultaneously quantales and covered semilattices. That is, they are quantales for which the multiplication equals the binary meet.

*Remark* 2.30. It may be useful to think about colax morphisms in the context of frames, though we will not do this. It is worth mentioning that a weak morphism of covered semilattices is automatically colax, so in particular a lax morphism of frames is automatically strong.

Note that, since every suplattice has a top element, which is necessarily a unit for $\wedge$, every frame is in fact a unital quantale. By Proposition 2.16, we can think of frames as exactly the idempotent affine unital quantales.

*Remark* 2.31. Most authors take morphisms of frames to be the strong morphisms as unital quantales, i.e. morphisms that preserve the top element. We will instead call these unital morphisms, in line with the terminology for quantales, as both the class of strong morphisms and the class of strong unital morphisms will be important to us, and we want to make the distinction clear.

**Theorem 2.32.** *Frames are reflective in quantales. That is, there is a functor* $\mathrm{Spec} : \mathbf{Qua} \to \mathbf{Frm}$, *and natural transformation* $\mathrm{id}_{\mathbf{Qua}} \to \mathrm{Spec}$, *with the following universal property: If* $\mathcal{Q}$ *is a quantale,* $\mathcal{F}$ *is a frame, and* $f : \mathcal{Q} \to \mathcal{F}$ *is a strong morphism of quantales, there exists a unique morphism of frames* $\mathrm{Spec}\, \mathcal{Q} \to \mathcal{F}$ *such that the following diagram commutes:*

$$\begin{array}{ccc} \mathcal{Q} & \longrightarrow & \mathrm{Spec}\, \mathcal{Q} \\ & \searrow_{f} & \downarrow \exists! \\ & & \mathcal{F} \end{array}$$

*Proof.* Rosenthal, p.44, Theorem 3.2.5.[14]. □

The proof of this theorem is nontrivial, though the commutative case is a good deal easier. Rosenthal proceeds by constructing three nuclei that encode algebraic properties of frames, and formally takes their join to obtain the minimal nucleus such that the corresponding quotient is a frame. Then, he shows that strong morphisms to frames lift across this nucleus—in effect, verifying that all such morphisms satisfy the lifting property of our Proposition 2.25.

We will suggestively refer to this nucleus on $\mathcal{Q}$ as $S \mapsto \sqrt{S}$.

*Remark* 2.33. The construction of $\sqrt{-}$ is abstract in general, but there are concrete ways to understand it in particular cases. For example, if every flat of $A$ is an ideal, then the flats $C$ with $\sqrt{C} = C$ are exactly the *semiprime ideals*, those satisfying $D \cdot D \subset C \implies D \subset C$ for any flat $D$. When $A$ is also commutative, finitary, and with unit, $\sqrt{C}$ is the radical of $C$ in the usual sense of commutative algebra.

For the reader who is not satisfied taking this construction as a black box, we give a detailed outline of its mechanics:

*Rosenthal's construction*

**Definition 2.34.** Let $\mathcal{Q}$ be a quantale, and $S \in \mathcal{Q}$. We say that $S$ is:

- **symmetric** if, whenever $U_1 \cdots U_n \leqslant S$, we also have $U_{\pi(1)} \cdots U_{\pi(n)} \leqslant S$ for any permutation $\pi$.

- **right-sided** if $S \cdot \top \leqslant S$, where $\top$ is the top element of $\mathcal{Q}$.

- **semiprime** if $T^2 \leqslant S$ implies $T \leqslant S$ for any $T \in \mathcal{Q}$.

- **localic** if S is symmetric, right-sided, and semiprime.

The goal is to demonstrate the minimum nucleus $j_{loc}$ on $\mathcal{Q}$ for which the quotient $\mathcal{Q}_{j_{loc}}$ is a frame. This turns out to be the same as finding the minimal nucleus on $\mathcal{Q}$ whose fixed points are all localic.

This proceeds in three parts. First, observe that an infimum of (symmetric / right-sided / semiprime) elements is (symmetric / right-sided / semiprime), so that we get three closure operators $j_c, j_r, j_e$. Next, verify that these closure operators are in fact nuclei, i.e. they are compatible with the multiplication of $\mathcal{Q}$.

The next step is the most non-constructive: define $j_{loc}$ as the join of nuclei $j_c \vee j_r \vee j_e$. Recall that nuclei are closed under pointwise infima, so the join $j_c \vee j_r \vee j_e$ is just the infimum of all nuclei $j$ with $j_c, j_r, j_e \leqslant j$. (For those who worry about empty infima, note that there is a maximum nucleus $j(a) = \top$.)

This shows that every quantale has a minimal frame quotient, and it provides some tools to construct it in certain situations. For example, Rosenthal shows that if $\mathcal{Q}$ is affine (two-sided), then an element is localic if and only if it is semiprime, and therefore $j_{loc}$ simply sends an element to the infimum of all semiprime elements above it. When applied to the quantale of ideals in a commutative ring, this recovers the usual notion of the radical.

These results are not quite enough to prove Theorem 2.32. The missing piece is to demonstrate that any morphism $f : \mathcal{Q} \to \mathcal{F}$ to a frame actually satisfies the lifting property of Proposition 2.25, i.e. is compatible with the nucleus $j_{loc}$. Rosenthal checks this separately for the nuclei $j_c, j_r, j_e$. The proof, while not particularly difficult, involves an application of the Adjoint Functor Theorem for lattices.

## 2.6 POINTLESS TOPOLOGY

The heart of pointless topology is a contravariant functor $\mathcal{O}$ on the category of topological spaces, sending a space $\mathcal{X}$ to its frame of opens $\mathcal{O}(\mathcal{X})$. This functor is full onto the category of frames and unital morphisms.

**Definition 2.35.** A frame is **spatial** if it is isomorphic to $\mathcal{O}(\mathcal{X})$ for some topological space $\mathcal{X}$.

The functor $\mathcal{O}$ has an adjoint pt, which we will describe in a moment, and this adjunction restricts to a duality between the category of spatial frames with unital morphisms, and the category of $T_0$ sober topological spaces.

*Remark* 2.36. To better resemble the theory of topological spaces, pointless topologists often work in the dual of the category of frames and unital morphisms, called the category of *locales*. In fact, the name "quantale" was originally a contraction of "quantum locale".[11]

*Remark* 2.37. The functor Spec gives us a way to construct frames, and frequently topological spaces, from quantales, and therefore from covered semigroups. However, to get topologically meaningful morphisms, we must take care to ensure that this functor only produces unital morphisms of frames. The most natural way to do this is to restrict to the category of covered monoids from the start, so that all morphisms will be unital. An alternative is to take the whole category of covered semigroups, but restrict to covering morphisms. These approaches differ in some situations, and coincide in others.

**Definition 2.38.** A **prime** or **point** of a frame $\mathcal{F}$ is an element $S \in \mathcal{F}$ such that $[\{S\}]$ is a prime ideal. Equivalently, we have $U \wedge V \leqslant S \implies U \leqslant S$ or $V \leqslant S$. The set of all points of $\mathcal{F}$ will be denoted by $\mathrm{pt}(\mathcal{F})$, and given the topology $\{D(S) \mid S \in \mathcal{F}\}$, where $D(S) = \{P \text{ prime} \mid U \not\leqslant P\}$.

*Remark* 2.39. Since all flats in a suplattice are principal with a unique generator, the prime elements of $\mathcal{F}$ are in bijection with the prime ideals (or prime filters) of $\mathcal{F}$. We will usually use "point" to refer to prime filters and ideals, and "prime" to refer to the prime elements of $\mathcal{F}$.

*Example* 2.40. If $\mathfrak{X}$ is a topological space, the primes of $\mathcal{O}(\mathfrak{X})$ are exactly the open sets whose complements are irreducible, and the space $\mathrm{pt}(\mathcal{O}(\mathfrak{X}))$ is the soberification of $\mathfrak{X}$.

**Definition 2.41.** A frame $\mathcal{F}$ has **enough points** if, for all $S \in \mathcal{F}$:

$$S = \bigwedge_{\substack{S \leqslant P, \\ P \text{ prime}}} P$$

Equivalently, whenever $S \not\leqslant T$, there is a prime $P$ such that $S \not\leqslant P$ and $T \leqslant P$.

**Proposition 2.42.** *A frame is spatial if and only if it has enough points.*

*Proof.* Johnstone, Section II.1.5 [6]. □

## 2.7 $\mathbb{F}_1$

**Definition 2.43.** $\mathbb{F}_1$ is defined to be the suplattice $\{0 \leqslant 1\}$, equipped with the usual multiplication $\cdot$. We regard $\mathbb{F}_1$ as a covered monoid (in fact, a unital quantale) via its suplattice closure, i.e. $[\emptyset] = [\{0\}] = \{0\}$, $[\{1\}] = [\mathbb{F}_1] = \mathbb{F}_1$.

$\mathbb{F}_1$ is often defined as simply the monoid $\{0, 1\}$, sometimes as the ordered monoid $\{0 \leqslant 1\}$. But it is our view that we should attach the correct closure operator to it, namely the one that comes from its suplattice structure. Specifically, we desire $[\emptyset] = \{0\}$, not $\emptyset$, which is what we would get by treating $\mathbb{F}_1$ as merely an ordered monoid.

$\mathbb{F}_1$ appears in many different guises:

- It is the frame of open sets of the one-point space.
- It is the quantale of flats of the trivial monoid.
- It is the initial object in the category of unital quantales.

One feature of $\mathbb{F}_1$ is that any set map $f : A \to \mathbb{F}_1$ is uniquely determined by either $f^{-1}(0)$ or $f^{-1}(1)$. This seemingly innocuous fact will help us to establish a useful dictionary between some familiar constructions.

**Definition 2.44.** If $A$ is a set and $f : A \to \mathbb{F}_1$ is a function, call $f$ **nontrivial** if $1 \in f(A)$.

If $A$ is a covered semigroup, we have the following correspondence of properties for a set map $f : A \to \mathbb{F}_1$:

| f | $K (= f^{-1}(0))$ | $F (= f^{-1}(1))$ |
|---|---|---|
| f is nontrivial | $K \neq A$ | $F \neq \emptyset$ |
| f is a closure map | $K = [K]$ | $a \in F, S \vdash a \implies F \cap S \neq \emptyset$ |
| f is lax | $x \in K$ or $y \in K \implies x \cdot y \in K$ | $x \cdot y \in F \implies x, y \in F$ |
| f is colax | $x \cdot y \in K \implies x \in K$ or $y \in K$ | $x, y \in F \implies x \cdot y \in F$ |

From this dictionary, we can see a few things immediately: The weak morphisms $A \to \mathbb{F}_1$ parametrize the flats of $A$, the lax morphisms $A \to \mathbb{F}_1$ parametrize the ideals of $A$, and the nontrivial strong morphisms $A \to \mathbb{F}_1$ parametrize the prime ideals of $A$. If $A$ is a semigroup with trivial coverage, the colax morphisms $A \to \mathbb{F}_1$ parametrize the subsemigroups of $A$, and if $A$ is a monoid with trivial coverage, then the unital colax morphism $A \to \mathbb{F}_1$ parametrize the submonoids of $A$. Finally, if $A$ is a meet-semilattice with the poset coverage, the colax morphisms $A \to \mathbb{F}_1$ parametrize the filters on $A$.

## 2.8 SPECTRA OF COVERED SEMIGROUPS

We can extend the functor Spec to the category of covered semigroups as follows. If $A$ is a covered semigroup, then $\mathcal{C}(A)$ is a quantale, and we define $\mathrm{Spec}\, A = \mathrm{Spec}\, \mathcal{C}(A)$. Since $\mathcal{Q} \to \mathcal{C}(\mathcal{Q})$ is an isomorphism for any quantale $\mathcal{Q}$, this is safe notation.

**Proposition 2.45.** *If $A$ is a covered semigroup, the prime ideals of $A$ are in bijection with the points of* Spec $A$.

*Proof.* This is a straightforward application of Theorem 2.32, Proposition 2.27, and the above dictionary. The prime ideals of $A$ are in bijection with the nontrivial strong morphisms of covered semigroups $A \to \mathbb{F}_1$, which are in bijection with the nontrivial strong morphisms of quantales $\mathcal{C}(A) \to \mathbb{F}_1$, which are in bijection with the nontrivial morphisms of frames $\mathrm{Spec}\, A \to \mathbb{F}_1$, which are in bijection with the primes of Spec $A$, which are in bijection with pt(Spec $A$). □

**Corollary 2.46.** *If $A$ is a covered semigroup,* Spec $A$ *is spatial if and only if, for every flat $C \in \mathcal{C}(A)$, we have:*

$$\sqrt{C} = \bigcap_{\substack{C \subset P, \\ P \text{ is a prime ideal}}} P$$

**Definition 2.47.** An element $S$ of a quantale $\mathcal{Q}$ is **compact** if, whenever $S \leqslant \bigvee \mathcal{S}$ for some collection $\mathcal{S} \subset \mathcal{Q}$, then $S \leqslant \bigvee \mathcal{S}'$ for some finite collection $\mathcal{S}' \subset \mathcal{S}$. A quantale is **algebraic** if every element is a join of compact elements.

**Theorem 2.48.** *Every algebraic quantale has a spatial spectrum.*

*Proof.* Rosenthal, p.59, Theorem 4.1.1. [14]. □

*Remark* 2.49. The proof of Theorem 2.48 requires the axiom of choice, and it is the only result in this work that does. If we are comfortable working with frames instead of topological spaces, the axiom of choice can be avoided entirely.

*Remark* 2.50. Not every $T_0$ sober space is the spectrum of an algebraic quantale, as such a spectrum necessarily has a basis of quasi-compact opens. However, we will see later, in the interlude, that all such topological spaces arise in this way.

**Definition 2.51.** A topological space is called **quasi-spectral** if it is $T_0$ and sober, and it has a basis of quasi-compact opens that is closed under intersection. A space is **spectral** if it is quasi-spectral and quasi-compact.

**Proposition 2.52.** *Let $A$ be a covered semigroup. If $A$ is finitary, then $\operatorname{Spec} A$ is spatial, and $\operatorname{pt}(\operatorname{Spec} A)$ is quasi-spectral. If $A$ has a unit, then $\operatorname{pt}(\operatorname{Spec} A)$ is spectral.*

*Proof.* Let $A$ be a finitary covered semigroup. Since every flat is a join of principal flats, Theorem 2.48 tells us that, to show $\operatorname{Spec} A$ is spatial, it is enough to show that every principal flat is compact. Suppose that $[\{a\}] \leqslant \bigvee \mathcal{S} = [\bigcup \mathcal{S}]$ for some collection $\mathcal{S}$ of flats. Since $A$ is finitary, there is a finite collection $\mathcal{S}' \subset \mathcal{S}$ such that $a \in [\bigcup \mathcal{S}'] = \bigvee \mathcal{S}'$, so $[\{a\}] \subset \bigvee \mathcal{S}'$, as desired.

We can identify $\{D(a) \mid a \in A\}$ with a basis of $\operatorname{pt}(\operatorname{Spec} A)$ of quasi-compact opens, and $D(a) \cap D(b) = D(a \cdot b)$, so it is closed under intersection, thus $\operatorname{pt}(\operatorname{Spec} A)$ is quasi-spectral. If $A$ has a unit, then $D(1)$ is the top element of $\operatorname{Spec} A$, therefore $\operatorname{pt}(\operatorname{Spec} A) = D(1)$ is quasi-compact, hence spectral. $\square$

# INTERLUDE: ALGEBRA OVER THE POWER MONAD

## N.1 HYPEROPERATIONS

**Definition N.1.** Let $A$ be a set, $n \geqslant 0$. An $n$-**ary hyperoperation** on $A$ is any set map $\star : A^n \to \mathcal{P}(A)$. (Note that we allow the range to include the empty set.)

**Definition N.2.** Let $(\mathcal{L}, \bigvee)$ be a suplattice, $n \geqslant 0$. An $n$-**ary lattice operation** on $\mathcal{L}$ is a map $\star : \mathcal{L}^n \to \mathcal{L}$ that is a morphism of suplattices $\mathcal{L} \to \mathcal{L}$ in each variable. Equivalently, the following should hold for all collections $\{S^i_{\alpha_i}\}$ of elements of $\mathcal{L}$:

$$\star\left(\bigvee_{\alpha_1} S^1_{\alpha_1}, \ldots, \bigvee_{\alpha_n} S^n_{\alpha_n}\right) = \bigvee_{\alpha_1, \ldots, \alpha_n} \star(S^1_{\alpha_1}, \ldots, S^n_{\alpha_n})$$

We will often work in the powerset rather than with $A$ directly. When it will not cause confusion, we will blur the distinction between the element $a \in A$ and the singleton $\{a\} \in \mathcal{P}(A)$.

Given an $n$-ary hyperoperation on $A$, there is a natural extension to $\mathcal{P}(A)$, given by:

$$\star(S_1, \ldots, S_n) = \bigcup_{s_i \in S_i} \star(s_1, \ldots, s_n)$$

**Proposition N.3.** *This extension induces a bijection between $n$-ary hyperoperations on $A$ and $n$-ary lattice operations on $(\mathcal{P}(A), \bigcup)$.*

*Proof.* We can check that any lattice operation on $(\mathcal{P}(A), \bigcup)$ is uniquely determined by its values on $A$ according to the above formula:

$$\star(S_1, \ldots, S_n) = \star(\bigcup_{s_1 \in S_1} \{s_1\}, \ldots, \bigcup_{s_n \in S_n} \{s_n\}) = \bigcup_{s_i \in S_i} \star(s_1, \ldots, s_n)$$

So we just need to verify that the formula always gives a lattice operation:

$$\star\left(\bigcup_{\alpha_1} S^1_{\alpha_1}, \ldots, \bigcup_{\alpha_n} S^n_{\alpha_n}\right) = \bigcup_{s_i \in \bigcup_{\alpha_i} S^i_{\alpha_i}} \star(s_1, \ldots s_n)$$

$$= \bigcup_{\alpha_1,...,\alpha_n} \bigcup_{s_i \in S^i_{\alpha_i}} \star(s_1, \ldots s_n) = \bigcup_{\alpha_1,...,\alpha_n} \star(S^1_{\alpha_1}, \ldots S^n_{\alpha_n})$$

□

*Remark* N.4. Note that we proved the $n = 1$ case much earlier, in Chapter 1.

*Example* N.5. A $\bigcup$-preserving map $\mathcal{P}(A) \to \mathcal{P}(A)$ is the same as a unary hyperoperation on $A$, i.e. a relation on $A$ or a directed graph with vertex set $A$.

*Remark* N.6. It is possible to develop this correspondence more abstractly using the tensor product of suplattices. As we have defined it, an $n$-ary lattice operation on $\mathcal{L}$ is exactly a morphism of suplattices $\mathcal{L}^{\otimes n} \to \mathcal{L}$. Then, Proposition N.3 follows from the fact that the functor **Set** → **Sup** sending $A$ to $\mathcal{P}(A)$ is a strong morphism of monoidal categories ($\mathcal{P}(A \times B) \cong \mathcal{P}(A) \otimes \mathcal{P}(B)$) and is left adjoint to the forgetful functor **Sup** → **Set**, as follows:

$$\mathrm{Hom}_{\mathbf{Set}}(A^n, \mathcal{P}(A)) \cong \mathrm{Hom}_{\mathbf{Sup}}(\mathcal{P}(A^n), \mathcal{P}(A)) \cong \mathrm{Hom}_{\mathbf{Sup}}((\mathcal{P}(A))^{\otimes n}, \mathcal{P}(A))$$

As a consequence of Proposition N.3, we can always study classical operations on lattices instead of hyperoperations. However, there is frequently some work to be done in relating properties of a hyperstructure to properties of its associated lattice. We will see shortly that associativity and commutativity of binary operations lift, as does a multiplicative hyperunit. However, later we will see that distributivity does not lift from semirings to lattice-enriched semirings.

A more basic counterexample is additive inverses: Although every element of $\mathbb{Z}$ has an additive inverse, the only elements of $\mathcal{P}(\mathbb{Z})$ with additive inverses are the singletons. In any case, the next section concerns itself with the algebraic properties that extend from $A$ to $\mathcal{P}(A)$ without incident.

### N.2 HYPERSEMIGROUPS AND HYPERMONOIDS

The definitions and results of this section appear to be folkloric. Some have been observed in online discussions by Trimble (for example, [16]).

**Definition N.7.** Let $\boxdot$ be a binary hyperoperation on a set $A$. We say that $\boxdot$ is **associative** if $a \boxdot (b \boxdot c) = (a \boxdot b) \boxdot c$ (as sets) for all $a, b, c \in A$. A set equipped with an associative binary hyperoperation is called a **hypersemigroup**.

**Proposition N.8.** *Let $A$ be a set, $\square$ any binary hyperoperation on $A$. Then $(A, \square)$ is a hypersemigroup if and only if $(\mathcal{P}(A), \bigcup, \square)$ is a quantale.*

*Proof.* One direction is easy: the associative law for the quantale $(\mathcal{P}(A), \bigcup, \square)$ implies the associative law for $(A, \square)$, a fortiori.

For the other direction, suppose that $\square$ is associative on $A$. By Proposition N.3, we know that $\square$ distributes over $\bigcup$, so all that remains to check is associativity. If $S, T, U$ are subsets of $A$, we have, by distributivity:

$$S \square (T \square U) = \left( \bigcup_{s \in S} s \right) \square \left( \bigcup_{t \in T, u \in U} t \square u \right)$$

$$= \bigcup_{s \in S, t \in T, u \in U} s \square (t \square u) = \bigcup_{s \in S, t \in T, u \in U} (s \square t) \square u$$

$$= \left( \bigcup_{s \in S, t \in T} s \square t \right) \square \left( \bigcup_{u \in U} u \right) = (S \square T) \square U$$

Since $\square$ is associative and distributes over $\bigcup$, we conclude that $(\mathcal{P}(A), \bigcup, \square)$ is a quantale. $\square$

**Definition N.9.** A **hypermonoid** is a tuple $(A, \square, 1)$, where $(A, \square)$ is a hypersemigroup, and $1$ is a nullary hyperoperation (i.e., a set $1 \in \mathcal{P}(A)$) such that $1 \square a = a \square 1 = a$ for all $a \in A$.

*Remark* N.10. There may be good reasons to generalize this by merely requiring $a \leqslant 1 \square a$ and $a \leqslant a \square 1$, though we will not explore this idea here.

**Proposition N.11.** *$(A, \square, 1)$ is a hypermonoid if and only if $(\mathcal{P}(A), \bigcup, \square, 1)$ is a unital quantale.*

*Proof.* If $1$ is a unit in $\mathcal{P}(A)$, then it is certainly a unit in $A$. Conversely, by Proposition N.8, all we need to check is that a unit in $A$ extends to a unit in $\mathcal{P}(A)$:

$$1 \square S = \bigcup_{s \in S} 1 \square s = \bigcup_{s \in S} s = S$$

Similarly, $S \square 1 = S$, so $1$ is a unit of $\mathcal{P}(A)$. $\square$

*Example* N.12. It is easy to check that every semigroup is a hypersemigroup, and every monoid is a hypermonoid.

*Example* N.13. If A is any set, then $(\mathcal{P}(A), \bigcup, \cap)$ is a frame, therefore a unital quantale, and so we have a hypermonoid structure on A given by:

$$1 = A$$

$$a \boxdot b = \{a\} \cap \{b\}$$

*Example* N.14. Similarly, if A is any set, we have a hypersemigroup structure on A, given by:

$$a \boxplus b = \{a\} \cup \{b\}$$

This is not, however, a hypermonoid: Since we are forced to define $\emptyset \boxplus S = \bigcup_{a \in \emptyset, b \in S} a \boxplus b = \emptyset$, there is no hyperunit when A has more than one element.

**Definition N.15.** A hypersemigroup or binary hyperoperation is **commutative** if $a \boxdot b = b \boxdot a$ for all $a, b \in A$.

**Proposition N.16.** *A binary hyperoperation $\boxdot$ is commutative on A if and only if it is commutative on $\mathcal{P}(A)$.*

*Proof.* If $\mathcal{P}(A)$ is commutative, then A is immediately commutative. If A is commutative, and $S, T \in \mathcal{P}(A)$, then:

$$S \boxdot T = \bigcup_{s \in S, t \in T} s \boxdot t = \bigcup_{t \in T, s \in S} t \boxdot s = T \boxdot S$$

$\square$

## N.3 COVERED HYPERSEMIGROUPS AND THEIR SPECTRA

Many basic results about covered semigroups extend immediately to hypersemigroups.

**Definition N.17.** Let $(A, \boxdot)$ be a hypersemigroup. A **coverage** on A is a closure such that $[S] \boxdot [T] \subset [S \boxdot T]$ for all $S, T \subset A$. A **covered hypersemigroup** is a hypersemigroup equipped with a coverage, and likewise a **covered hypermonoid** is a hypermonoid equipped with a coverage.

Since a coverage on a hypersemigroup $(A, \boxdot)$ is exactly a nucleus on the quantale $(\mathcal{P}(A), \bigcup, \boxdot)$, the lattice of flats $\mathcal{C}(A)$ is naturally a quantale, allowing us to extend the functor Spec to covered hypersemigroups in the same way we did for covered semigroups.

**Proposition N.18.** *The spectrum of a finitary covered hypersemigroup is spatial.*

*Proof.* The proof is identical to the first half of the proof of the similar proposition for covered semigroups in Chapter 2. □

**Proposition N.19.** *If $\mathfrak{X}$ is a $T_0$, sober topological space with a basis of quasi-compact opens, then $\mathfrak{X}$ is the spectrum of a finitary covered hypersemigroup.*

*Proof.* Let A be the set of all quasi-compact opens of $\mathfrak{X}$. For $a, b \in A$, define $a \boxdot b = \{x \in A \mid x \subset a \cap b\}$, and give A the closure defined by $S \dashv a$ if $a \subset \bigcup S$. It is not difficult to check that A satisfies the axioms for a finitary covered hypersemigroup.

Furthermore, if C is a flat and $x \subset \bigcup C$, then $x \in C$, so $C \mapsto \bigcup C$ is an isomorphism of suplattices $\mathcal{C}(A) \to \mathcal{O}(\mathfrak{X})$, which preserves the multiplication and is therefore an isomorphism of quantales. It follows that $\mathcal{C}(A)$ is already a frame, so $\operatorname{Spec} A = \mathcal{C}(A) \cong \mathcal{O}(\mathfrak{X})$. □

### n.4 subdistributivity

One basic algebraic property that almost never lifts to the powerset is distributivity. For example:

$$\{0, 1\} \cdot (\{1\} + \{1\}) \subsetneq \{0, 1\} \cdot \{1\} + \{0, 1\} \cdot \{1\}$$

However, there is a much better-behaved property:

**Definition N.20.** If $\boxdot, \boxplus$ are two binary hyperoperations on a set, we say that $\boxdot$ **subdistributes** over $\boxplus$ if the following two conditions hold:

$$a \boxdot (b \boxplus c) \subset (a \boxdot b) \boxplus (a \boxdot c)$$
$$(a \boxplus b) \boxdot c \subset (a \boxdot c) \boxplus (b \boxdot c)$$

If $\boxdot, \boxplus$ are two binary lattice operations, we say that $\boxdot$ **subdistributes** over $\boxplus$ if the following two conditions hold:

$$S \boxdot (T \boxplus U) \leqslant (S \boxdot T) \boxplus (S \boxdot U)$$
$$(S \boxplus T) \boxdot U \leqslant (S \boxdot U) \boxplus (T \boxdot U)$$

**Proposition N.21.** *If $\boxdot$ subdistributes over $\boxplus$ in A, then $\boxdot$ subdistributes over $\boxplus$ in $\mathcal{P}(A)$.*

*Proof.*
$$S \odot (T \boxplus U) = \bigcup_{s,t,u} s \odot (t \boxplus u) \subset \bigcup_{s,t,u} (s \odot t \boxplus s \odot u)$$
$$\subset \bigcup_{s,s',t,u} (s \odot t \boxplus s' \odot u) = S \odot T \boxplus S \odot U$$

The other condition holds by symmetry. □

*Remark* N.22. It is no coincidence that subdistributivity is equivalent to the condition that multiplication by a fixed element is a colax morphism with respect to $\boxplus$. The next chapter will examine this more closely.

# 3

# QUANTIC SEMIRINGS AND HYPERSEMIRINGS

## 3.1 QUANTIC SEMIRINGS

**Definition 3.1.** A **quantic semiring** is a tuple $(\mathcal{S}, \bigvee, \cdot, +, 0)$ such that:

- $(\mathcal{S}, \bigvee, \cdot)$ is a quantale.
- $(\mathcal{S}, +, 0)$ is a commutative monoid.
- $+$ distributes over non-empty countable joins.
- For all $S, T, U \in \mathcal{S}$, $S \cdot (T + U) \leqslant S \cdot T + S \cdot U$ and $(S + T) \cdot U \leqslant S \cdot U + T \cdot U$.
- For all $S \in \mathcal{S}$, $0 \cdot S \leqslant 0$ and $S \cdot 0 \leqslant 0$.

*Remark* 3.2. The last axiom may be surprising, but it cannot easily be strengthened. Indeed, in the powerset of a ring, we have $\emptyset \cdot \{0\} = \emptyset$, not $\{0\}$. We can see this axiom as a form of sub-distributivity that applies to the empty sum.

*Remark* 3.3. In practice, the third axiom may be strengthened, so that $+$ distributes over all non-empty joins. But it is too much to assume that it distributes over the empty join, as this fails for one of our principal examples: quantales.

**Definition 3.4.** If $\mathcal{S}, \mathcal{S}'$ are quantic semirings, then a morphism of quantic semirings $f : \mathcal{S} \to \mathcal{S}'$ is a function on the underlying sets such that:

- For all $\mathcal{T} \subset \mathcal{S}$, $f(\bigvee \mathcal{T}) = \bigvee f(\mathcal{T})$.
  ($f$ is a morphism of suplattices.)

- For all $S, T \in \mathcal{S}$, $f(S) \cdot f(T) = f(S \cdot T)$.
  ($f$ is a strong morphism of multiplicative quantales.)

- For all $S, T \in \mathcal{S}$, $f(S + T) \leqslant f(S) + f(T)$, and $f(0) = 0$.
  ($f$ is a colax morphism of additive ordered monoids.)

As with $\wedge$, we use the lower-case $\vee$ to denote the binary join.

**Proposition 3.5.** *Let $(\mathcal{Q}, \bigvee, \cdot)$ be a quantale with bottom element $0$. Then $(\mathcal{Q}, \bigvee, \cdot, \vee, 0)$ is a quantic semiring. Furthermore, a morphism of quantales is a morphism of quantic semirings.*

*Proof.* It is clear that $(\mathcal{Q}, \vee, 0)$ is a commutative monoid. $\vee$ also distributes over non-empty joins: $S \vee \bigvee_{T \in \mathcal{T}} T = \bigvee_{T \in \mathcal{T}} S \vee \bigvee_{T \in \mathcal{T}} T = \bigvee_{T \in \mathcal{T}} (S \vee T)$. Furthermore, $0 \cdot S = 0$ and $S \cdot 0 = 0$ for all $S$, and the fact that $\cdot$ distributes over $\bigvee$ means that it also distributes over $\vee$.

To see the statement about morphisms, note that a morphism of suplattices preserves joins, and is therefore strong with respect to $\vee$. □

The above proposition allows us to interpret quantales as quantic semirings, which we will do from now on.

**Definition 3.6.** An element $M \in \mathcal{S}$ is **monoidal** if $0 \leqslant M$ and $M + M = M$. The collection of monoidal elements will be denoted $\mathcal{M}(\mathcal{S})$.

*Remark* 3.7. The sum of two monoidal elements is monoidal, but the product of two monoidal elements generally isn't.

**Proposition 3.8.** *A quantic semiring $\mathcal{S}$ is a quantale if and only if every element of $\mathcal{S}$ is monoidal.*

*Proof.* This is just the dual of the multiplicative version in Chapter 2. Every element of a quantale is clearly monoidal, and conversely if the conditions are satisfied, we have, for all $S, T \in \mathcal{S}$:

$$S + T \leqslant (S \vee T) + (S \vee T) = S \vee T$$
$$S = S + 0 \leqslant S + T$$
$$T = 0 + T \leqslant S + T$$

It follows that $S + T = S \vee T$, so $\mathcal{S}$ is a quantale. □

Next, we will give a universal way to "monoidify" a quantic semiring.

## 3.2 THE QUANTIC SPECTRUM

**Definition 3.9.** For $S \in \mathcal{S}$, define $0 * S = 0$ and $n * S = \underbrace{S + \cdots + S}_{n \text{ times}}$ for $n \geqslant 1$.

We also define:

$$[S] = \bigvee_{n \geqslant 0} n * S$$

*Remark* 3.10. It may be useful to the reader to note the similarities with the Kleene star of computer science. In fact, if A is a ring or monoid and $\mathcal{S} = \mathcal{P}(A)$, then [S] is exactly the multiplicative monoid generated by S.

**Proposition 3.11.** $S \mapsto [S]$ *is a nucleus on* $(\mathcal{S}, \bigvee, \cdot)$, *whose fixed points are exactly the monoidal elements.*

*Proof.* The first statement follows from sub-distributivity:

$$[S] \cdot [T] = \bigvee_{n \geqslant 0} n * S \cdot \bigvee_{m \geqslant 0} m * T$$

$$= \bigvee_{n,m \geqslant 0} (n * S) \cdot (m * T) \leqslant \bigvee_{n,m \geqslant 0} (nm) * (S \cdot T)$$

$$= \bigvee_{k \geqslant 0} k * (S \cdot T) = [S \cdot T]$$

Clearly, if S is monoidal then $[S] = S$. The fact that $[S]$ is always monoidal follows from $0 \leqslant [S]$ and the fact that $+$ distributes over non-empty countable joins:

$$[S] + [S] = \bigvee_{n \geqslant 0} n * S + \bigvee_{m \geqslant 0} m * S$$

$$\bigvee_{n,m \geqslant 0} (n * S + m * S) = \bigvee_{n,m \geqslant 0} (n + m) * S$$

$$= \bigvee_{k \geqslant 0} k * S = [S]$$

□

**Proposition 3.12.** $[-] : \mathcal{S} \to \mathcal{M}(\mathcal{S})$ *is a morphism of quantic semirings.*

*Proof.* Since $[-]$ is a nucleus, this is a strong morphism of quantales, so all we need to check is that $[S] + [T] \geqslant [S + T]$:

$$[S] + [T] = \bigvee_{n \geqslant 0} n * S + \bigvee_{m \geqslant 0} m * T$$

$$= \bigvee_{n,m \geqslant 0} (n * S + m * T) \geqslant \bigvee_{k \geqslant 0} (k * S + k * T)$$

$$= \bigvee_{k \geqslant 0} k * (S + T) = [S + T]$$

□

**Proposition 3.13.** *Let $\mathcal{S}$ be a quantic semiring, $\mathcal{Q}$ a quantale, and $f : \mathcal{S} \to \mathcal{Q}$ a morphism of quantic semirings. There is a unique strong morphism of quantales $[f] : \mathcal{M}(\mathcal{S}) \to \mathcal{Q}$ making the following diagram commute:*

$$\begin{array}{ccc} \mathcal{S} & \longrightarrow & \mathcal{M}(\mathcal{S}) \\ & f \searrow & \downarrow [f] \\ & & \mathcal{Q} \end{array}$$

*Proof.* By the results of Chapter 2, it is enough to show that $f(S) = f([S])$ for all $S \in \mathcal{S}$. Since $f(S) \leqslant f([S])$, it suffices to show that $f(n * S) \leqslant f(S)$ for all $n$:

- $f(0 * S) = f(0) = 0 \leqslant f(S)$.

- $f(n * S) = f(\underbrace{S + \cdots + S}_{n \text{ times}}) \leqslant \underbrace{f(S) + \cdots + f(S)}_{n \text{ times}} = f(S)$.

□

This shows that quantales are a reflective subcategory of quantic semirings: to every quantic semiring, there is a universal associated quantale.

**Definition 3.14.** For a quantic semiring $\mathcal{S}$, we call the quantale $\mathcal{M}(\mathcal{S})$ the **quantic spectrum** of $\mathcal{S}$. We define the classical spectrum of $\mathcal{S}$ by composition with $\mathcal{M}$: $\operatorname{Spec} \mathcal{S} = \operatorname{Spec} \mathcal{M}(\mathcal{S})$.

3.3 HYPERSEMIRINGS

**Definition 3.15.** A **hypersemiring** is a tuple $(A, \boxdot, \boxplus, 0)$ such that:

- $(A, \boxdot)$ is a hypersemigroup.

- $(A, \boxplus, 0)$ is a commutative hypermonoid.

- For all $a, b, c \in A$, $a \boxdot (b \boxplus c) \subset (a \boxdot b) \boxplus (a \boxdot c)$ and $(a \boxplus b) \boxdot c \subset (a \boxdot c) \boxplus (b \boxdot c)$.

- For all $a \in A$, $0 \boxdot a \subset 0$ and $a \boxdot 0 \subset 0$.

**Proposition 3.16.** *Let $A$ be a set, $\boxdot$ and $\boxplus$ binary hyperoperations and $A$, and $0 \subset A$. Then $(A, \boxdot, \boxplus, 0)$ is a hypersemiring if and only if $(\mathcal{P}(A), \bigcup, \boxdot, \boxplus, 0)$ is a quantic semiring.*

*Proof.* The last two conditions follow a fortiori, and the first condition follows from the results of the interlude. Note that $⊞$ is a lattice operation by the interlude, so it necessarily distributes over all joins. Then, if $(\mathcal{P}(A), \bigcup, □, ⊞, 0)$ is a quantic semiring, $(\mathcal{P}(A), \bigcup, +, 0)$ is a commutative unital quantale, and, by the results of the interlude, the second condition follows.

Conversely, if $(A, □, ⊞, 0)$ is a hypersemiring, then $(\mathcal{P}(A), \bigcup, □)$ is a quantale by the interlude, $(\mathcal{P}(A), \bigcup, ⊞, 0)$ is a commutative unital quantale by the interlude, sub-distributivity follows from the interlude, and the final condition is easily shown:

$$0 □ S = \bigcup_{s \in S} 0 □ s \subset \bigcup_{s \in S} 0 = 0$$

Similarly, $S □ 0 \subset 0$. □

**Definition 3.17.** A morphism of hypersemirings is a function $f : A \to B$ that lifts to a morphism $\mathcal{P}(A) \to \mathcal{P}(B)$ of quantic semirings. Equivalently, it is a strong morphism of multiplicative hypersemigroups, and a colax morphism of additive hypermonoids.

As a result of Proposition 3.16, we can assign a canonical coverage to a hypersemiring $A$, namely the monoid nucleus $S \mapsto [S]$ on the quantic semiring $\mathcal{P}(A)$. Since every hypersemiring is then naturally a covered hypersemigroup, it makes sense to talk about its spectrum.

**Definition 3.18.** An **ideal** of a hypersemiring $A$ is a subset $I \subset A$ satisfying:

- $0 \subset I$.
- $I + I = I$.
- $a \in I$ or $b \in I \implies a \cdot b \subset I$.

An ideal $I$ is **prime** if it also satisfies:

- $I \neq A$.
- $a \cdot b \subset I \implies a \in I$ or $b \in I$.

**Theorem 3.19.** *The spectrum of a hypersemiring is spatial. Its points are given by the prime ideals of $A$, and it has a basis of quasi-compact opens given by the sets $D(a) = \{I \mid a \notin I\}$ for $a \in A$.*

*Proof.* To see that $\operatorname{Spec} A$ is spatial, it is enough to check that the coverage on $A$ is finitary. Suppose that $a \in [S] = \bigvee_{n \geqslant 0} n * S$. Then $a \in n * S$ for some $n$. Either $n = 0$, in which case $a \in 0 = [\emptyset]$, or $n > 0$, in which case there are $s_1, \ldots, s_n$ such that $a = s_1 + \cdots + s_n$. Then $a \in [\{s_1, \ldots, s_n\}]$.

To compute the points, note that we have $0 \subset I$ and $I + I = I$ if and only if $I$ is a flat of the coverage on $A$. Thus the prime ideals of $A$ are exactly the prime ideals of $\mathcal{M}(\mathcal{P}(A))$, which are the points of $\operatorname{Spec} \mathcal{M}(\mathcal{P}(A)) = \operatorname{Spec} A$.

$\square$

# 4

# EXAMPLES AND APPLICATIONS

## 4.1 THE SPECTRUM OF A KRASNER HYPERRING

The primary example motivating our development of quantic semirings was that of the spectrum of a Krasner hyperring. Specifically, we wanted to show that existing definitions could be arrived at through a process more canonical than imitating the definitions from standard commutative ring theory.

To show that we have been successful, we will show that our approach is consistent with the standard definitions. Our main reference is Jun's recent paper [7], though these ideas go back further (e.g. [13]).

Let $A$ be a Krasner hyperring—that is, a unital hypersemiring with a commutative, single-valued multiplication, satisfying two axioms for inverses:

- For every $a \in A$, there is a unique $b \in A$, denoted $-a$, such that $0 \in a \boxplus b$.

- For every $a, b, c \in A$, $a \in b \boxplus c \iff c \in a \boxplus (-b)$.

*Remark* 4.1. The above two axioms say that $(A, \boxplus, 0)$ is a *canonical hypergroup*. We will not actually use these two axioms, but they appear to be important for other theoretical concerns, and we include them for completeness.

*Remark* 4.2. Every hypersemiring with single-valued multiplication is automatically strongly distributive, so in particular every Krasner hyperring is strongly distributive. This is ususally taken as an axiom for Krasner hyperrings, but it follows in our case from sub-distributivity.

*Remark* 4.3. We have omitted the axiom that the addition be non-empty, but this follows from the existence of additive inverses.

Recall that we defined the coverage on a hypersemiring as follows: The flats are subsets $C \subset A$ such that $0 \in C$ and $C + C = C$. An ideal is then a flat $I$ such that $A \cdot I = I$.

In fact, this is equivalent to other definitions in the literature of an ideal of a Krasner hyperring. For example, the definition in [7] is that $I \neq \emptyset$ and $x, y \in I, a \in A \implies x - ay \subset I$. But $x - ay = x + a(-y) = x + a(-1)y$, since the multiplication is single-valued, so this precisely asserts that I is a flat (in our sense) and closed under multiplication by A.

It is then straightforward that our notion of prime ideal is equivalent to the usual one. We can also compare the topologies on pt(Spec A): [7] gives them a basis consisting of the principal opens $D(a)$, while our machinery does precisely the same thing. (See, for example, Proposition 2.52 and its proof.)

One notable insight is that the spectrum of a Krasner hyperring is spectral precisely because of two facts: first, the multiplication is single-valued, which guarantees that the standard basis is closed under intersection, and second, there is a unit for the multiplication, which guarantees that the spectrum itself belongs to the standard basis. This puts Krasner hyperrings among a broader class of examples: any hypersemiring such that the product of any two elements is a finite set must have a quasi-spectral spectrum, and if it furthermore has a finite hyper-unit, its spectrum is a spectral space.

## 4.2 NUCLEI ON ORDERED BLUEPRINTS

The ordered blueprints of Lorscheid, as defined in [10], provide a rich set of examples of our theory. Unlike in the case of commutative rings, where the various nuclei all give the same notion of ideal, on an ordered blueprint there may be several reasonable choices for a nucleus. Below, we explore three that we believe to be of interest, and discuss the implications of combining them.

Briefly, an **ordered blueprint** $(A, \cdot, \leqslant)$ is a multiplicative commutative monoid with absorbing element 0, together with a preorder $\leqslant$ on the semiring $\mathbb{N}[A]^+$ of formal finite sums of elements of A, that is compatible with the addition and multiplication on $\mathbb{N}[A]^+$, and also satisfies the following two axioms:

- If $\sum$ denotes the empty sum, $0 \leqslant \sum$ and $\sum \leqslant 0$.

- If $a \leqslant b$ and $b \leqslant a$ for $a, b \in A$, then $a = b$.

The relation $\leqslant$ is called a **preaddition**. If $\sum_i a_i \leqslant \sum_j b_j$ and $\sum_j b_j \leqslant \sum_i a_i$, we write $\sum_i a_i \equiv \sum_j b_j$. A blueprint is called **algebraic** if $\leqslant$ coincides with $\equiv$, i.e. if $\leqslant$ is symmetric.

Throughout the section, A and B will be ordered blueprints. Recall that a morphism of ordered blueprints is a morphism of multiplicative monoids $f : A \to B$ such that $\sum_i a_i \leqslant \sum_j b_j$ implies $\sum_i f(a_i) \leqslant \sum_j f(b_j)$.

*The algebraic nucleus*

**Proposition 4.4.** *The following relation gives a coverage on* $A$: $S \vdash a$ *if there are* $s_i, s'_j \in S$ *such that* $a + \sum_i s_i \equiv \sum_j s'_j$.

*Proof.* First, we check that $\vdash$ is a closure relation. Since $a \equiv a$, we certainly have $S \vdash a$ for $a \in S$. Suppose that $T \vdash s$ for every $s \in S$, and $S \vdash a$. Then there are $s_i, s'_j \in S$, $t_{ik}, t'_{il}, u_{jm}, u'_{jn} \in T$ such that:

$$a + \sum_i s_i \equiv \sum_j s'_j$$

$$s_i + \sum_k t_{ik} \equiv \sum_l t'_{il}$$

$$s'_j + \sum_m u_{jm} \equiv \sum_n u'_{jn}$$

Combining these gives us:

$$a + \sum_{il} t'_{il} + \sum_{jm} u_{jm}$$

$$\equiv a + \sum_i s_i + \sum_{ik} t_{ik} + \sum_{jm} u_{jm}$$

$$\equiv \sum_j s'_j + \sum_{ik} t_{ik} + \sum_{jm} u_{jm}$$

$$\equiv \sum_{jn} u'_{jn} + \sum_{ik} t_{ik}$$

So $T \vdash a$, and $\vdash$ is a closure relation.

To check that we have a coverage, we just need to show that $S \vdash a \implies x \cdot S \vdash x \cdot a$. But this follows immediately from multiplicativity of $\equiv$. $\square$

**Definition 4.5.** The above coverage is called the **algebraic coverage** or **algebraic nucleus**. Its flats and ideals will be called **algebraic flats** and **algebraic ideals** respectively.

Given a morphism $f : A \to B$ of ordered blueprints, it is not hard to see that $f$ is a morphism of covered monoids when we given both $A$ and $B$ the algebraic coverage. Indeed, if $a + \sum_i s_i \leq \sum_j s'_j$ and $a + \sum_i s_i \geq \sum_j s'_j$, then $f(a) + \sum_i f(s_i) \leq \sum_j f(s'_j)$ and $f(a) + \sum_i f(s_i) \geq \sum_j f(s'_j)$.

**Proposition 4.6.** *The ideals of the algebraic nucleus are exactly the k-ideals in the sense of Lorscheid, and the spectrum of A as a monoid with this coverage coincides with the prime spectrum as an ordered blueprint.*

*Proof.* A k-ideal is defined to be a subset $I \subset A$ with $0 \in I$, $I \cdot A = I$, such that if $c + \sum_i a_i \equiv \sum_j b_j$ with $a_i, b_j \in I$, then $c \in I$. The last condition is equivalent to saying that $I$ is an algebraic flat. Together with the condition $I \cdot A$, this says exactly that $I$ is an algebraic ideal. So, to compare k-ideals to algebraic ideals, it suffices to show that $0 \in I$ whenever $I$ is an ideal of the algebraic nucleus. But it follows from the definition that $\emptyset \vdash 0$, so every algebraic flat and in particular every algebraic ideal contains $0$. To compare the spectra, it is enough to notice that the points are the same (prime ideals) and the topologies are the same (generated by the principal sets $D(a)$). □

*Remark 4.7.* As we see, the algebraic nucleus is closely related to Lorscheid's definintion of the prime spectrum of an ordered blueprint. However, since it depends only on the algebraic core of A—the symmetric relations only— it is likely to be the wrong definition in many cases where the ordering is critical. In particular, if we take a Krasner hyperring qua ordered blueprint, the algebraic nucleus has too many ideals and prime ideals. We will fix this problem in the next example.

*The monomial nucleus*

**Proposition 4.8.** *The following relation gives a coverage on A: $S \vdash a$ if there are $s_i \in S$ such that $a \leqslant \sum_i s_i$.*

*Proof.* The proof is similar to the previous example, but simpler: If $a \leqslant \sum_i s_i$, and $s_i \leqslant \sum_j t_{ij}$, then $a \leqslant \sum_{ij} t_{ij}$. □

**Definition 4.9.** The above coverage is called the **monomial coverage** or **monomial nucleus**, its flats and ideals the **monomial flats** and **monomial ideals**.

As with the algebraic coverage, morphisms of ordered blueprints respect the monomial coverage, since $a \leqslant \sum_i s_i \implies f(a) \leqslant \sum_i f(s_i)$.

**Proposition 4.10.** *If A is a Krasner hyperring, then the ideals of A as a hyperring coincide with the monomial ideals of A as an ordered blueprint.*

*Proof.* It is enough to show that the monomial coverage on A coincides with the canonical coverage constructed in Chapter 3. Recall that the sub-addition on A is defined to be the one generated by relations of the form $a \leqslant b + c$, whenever $a \in b \boxplus c$. By induction, we can see that any subaddition containing these relations must contain all relations of the form $\sum_i a_i \leqslant \sum_j b_j$, whenever $\boxplus_i a_i \subset \boxplus_j b_j$. Furthermore, the collection of all such relations forms a subaddition, and is therefore the subaddition associated to A.

From there, the result is almost immediate: the monomial relations are those of the form $a \leqslant \sum_i s_i$, where $a \in \sum_i s_i$, and these relations give exactly the closure relation of Chapter 3. □

*Remark* 4.11. As we can see from the example of hyperrings, the monomial nucleus is useful in practice. However, it has the shortcoming that it depends only on the monomial relations of A, and may therefore be too crude in situations where more complicated relations are involved. In particular, for algebraic blueprints the monomial nucleus is more crude than the algebraic nucleus.

*The order nucleus*

**Proposition 4.12.** *The following relation gives a coverage on A: $S \vdash a$ if there is $s \in S$ such that $a \leqslant s$.*

*Proof.* This is a closure operator by our discussion of posets in Chapter 1. Since ordered blueprints are in particular ordered monoids, this is compatible with the multiplication as well. □

**Definition 4.13.** The above coverage is called the **order coverage** or **order nucleus**, its flats and ideals (recall Definition 2.5) the **order flats** and **order ideals**. If $S \subset A$, we define $S^{\leqslant} = \{a \leqslant s \mid s \in S\}$ to be the flat generated by S.

The spirit of this example is different from the other two. We are ignoring the additive structure on A, taking only its structure as an ordered monoid, so it would seem that we cannot extract much useful information. However, any nucleus that respects the order structure of A (for example, the monomial nucleus) will be a specialization of the order nucleus, so it is worth keeping in mind as an intermediate step. There is an analogy in Lorscheid's work: he uses the spectrum of an algebraic blueprint as a monoid to define the "subcanonical spectrum"; we propose that it may be interesting to consider the spectrum of an ordered blueprint as an ordered monoid.

There is a particular case in which the order nucleus provides a strong connection to our work in Chapter 3.

**Definition 4.14.** An ordered blueprint A is **associative** if, whenever $a \leqslant s + t + u$, there is some $b \leqslant s + t$ such that $a \leqslant b + u$.

*Remark* 4.15. We use the word "associative" for this property because it is actually (in some sense) equivalent to the associativity of addition for Krasner hyperrings: if A is a Krasner hyperring qua ordered blueprint, then $a \leqslant s + t + u$ implies $a \in s \boxplus t \boxplus u$ (we showed this in the proof of Proposition 4.10). So $a \in (s \boxplus t) \boxplus u = \bigcup_{b \in s \boxplus t} b \boxplus u$, so A is associative.

On the other hand, we will see in the proof below that associativity of a literal addition follows directly from this property.

**Proposition 4.16.** *If A is an associative ordered blueprint, the following operations make the order flats of A into a quantic semiring:*

$$\bigvee_\alpha S_\alpha = \bigcup_\alpha S_\alpha$$

$$S \odot T = (S \cdot T)^{\leqslant}$$

$$S \oplus T = \{a \mid a \leqslant s + t \text{ for some } s \in S, t \in T\}$$

$$0 = \{a \mid a \leqslant \emptyset\}$$

*Proof.* First, it is straightforward to check that these operations all produce downward-closed sets. Furthermore, since $S \mapsto S^{\leqslant}$ is a nucleus, it is immediate that the downward-closed sets form a multiplicative quantale. We will check the other properties directly:

- 0 is an additive identity: if $a \leqslant s + t$ with $t \in 0$, then $a \leqslant s + t \leqslant s + \emptyset = s$, and $a \in S$. It follows that $S \oplus 0 \subset S$. But since there is some $z \in A$ with $z \equiv \sum_\emptyset$, we have $a \leqslant s \implies a \leqslant s + z$, so $S \subset S \oplus 0$.

- $\oplus$ distributes over joins: We have $\bigcup_\alpha S_\alpha \oplus T = \{a \mid a \leqslant s + t \text{ for some } s \in \bigcup_\alpha S_\alpha, t \in T\} = \bigcup_\alpha \{a \mid a \leqslant s + t \text{ for some } s \in S_\alpha, t \in T\} = \bigcup_\alpha (S_\alpha \oplus T)$.

- Sub-distributivity: If $a \in S \odot (T \oplus U)$, then there are $s \in S, t \in T, u \in U$ such that $x \leqslant t + u$ and $a \leqslant s \cdot x$. But then $a \leqslant s \cdot t + s \cdot u$, so $a \in (S \odot U) \oplus (T \odot U)$.

- $S \odot 0 \subset 0$: If $t \leqslant \emptyset$, then $s \cdot t \leqslant \emptyset$, so $S \odot 0 \subset 0$.

- $\oplus$ is associative: if $a \in S \oplus (T \oplus U)$, then there are $s \in S, t \in T, u \in U$ with $a \leqslant s + t + u$. By associativity, there is some $x \in A$ with $x \leqslant s + t$ and $a \leqslant x + u$, so $a \in (S \oplus T) \oplus U$, therefore $S \oplus (T \oplus U) \subset (S \oplus T) \oplus U$. Similarly, $(S \oplus T) \oplus U \subset S \oplus (T \oplus U)$.

**Proposition 4.17.** *If* $f : A \to B$ *is a morphism of associative ordered blueprints, then* $F(S) = f(S)^{\leqslant}$ *defines a morphism of the associated quantic semirings.*

*Proof.* We check:

- F preserves joins:

$$F(\bigcup_\alpha S_\alpha) = f(\bigcup_\alpha S_\alpha)^{\leqslant} = (\bigcup_\alpha f(S_\alpha))^{\leqslant}$$

$$= \bigcup_\alpha f(S_\alpha)^{\leqslant} = \bigcup_\alpha F(S_\alpha)$$

- F is multiplicative: $F(S \odot T) = F(S \cdot T) = f(S \cdot T)^{\leqslant} = (f(S) \cdot f(T))^{\leqslant} = f(S) \odot f(T)$.

- F is colax-additive: Since f is a morphism of ordered blueprints, $a \leqslant s + t \implies f(a) \leqslant f(s) + f(t)$, and therefore $F(S \oplus T) = f(S \oplus T)^{\leqslant} \subset f(S) \oplus f(T) \subset F(S) \oplus F(T)$.

□

**Proposition 4.18.** *If* A *is an associative ordered blueprint, its monomial nucleus is equal to the composition of the order nucleus with the nucleus* M *on the quantic semiring of order flats.*

*Proof.* This follows from the fact that any relation $a \leqslant \sum_i s_i$ in an associative ordered blueprint is generated by relations of the form $a \leqslant b + c$. □

*Remark 4.19.* One interesting thing about this construction is that if the ordering on A is discrete, every subset of A is an order flat. In particular, if A is an algebraic blueprint or a Krasner hyperring, we are just passing to the powerset.

*Remarks on composite nuclei*

It is a useful fact that the collection of all nuclei on a semigroup is a suplattice, as it means that we can combine arbitrary nuclei in a universal manner. This combination may have nontrivial properties, since it is constructed in an abstract manner as the pointwise infimum of upper bounds.

So, for example, there is a canonical minimal nucleus sitting over both the algebraic nucleus and the order nucleus, and one (most likely different) sitting over the algebraic and monomial nuclei. These are good candidates for defining "order ideals" for ordered blueprints that would share features of the ideal theory of algebraic blueprints (k-ideals) and the ideal theory of hyperrings.

There do not appear to be one-line descriptions of these nuclei, but they can most likely be described in terms of finite chains of relations, alternating between algebraic relations and monomial/order relations.

## 4.3 THOUGHTS AND FURTHER DIRECTIONS

*Polynomials over hyperrings: double distributivity*

If $A$ is a hypersemiring, then there is a canonical multiplication and addition on the set $\mathcal{P}(A)[X]$ of power-polynomials, polynomials with coefficients in $\mathcal{P}(A)$: simply define $(\sum_i a_i X^i) \boxplus (\sum_i b_i X^i) = \sum_i (a_i \boxplus b_i) X^i$ and $(\sum_i a_i X^i) \boxdot (\sum_j b_j X^j) = \sum_n \sum_{i+j=n} (a_i \boxdot b_j) X^n$. This pulls back to a hyperstructure on $A[X]$: the hyperproduct of two polynomials is just the set of all polynomials contained in the power-polynomial product, and similarly with addition.

Unfortunately, $A[X]$ is not, in general, a hypersemiring, because the multiplication may easily fail to be associative. And there is another problem: the evaluation map $e_t : A[X] \to A$ sending $p(X)$ to $p(t)$ may fail to be a multiplicative homomorphism in a meaningful sense.

However, both of these problems go away if we impose the condition of *double distributivity*, which says that $(a+b)(c+d) = ac + ad + bc + bd$ for all $a, b, c, d \in A$. By induction, it implies that $(\sum_i a_i)(\sum_j b_j) = \sum_{ij} a_i b_j$, so it is in some sense the strongest possible form of distributivity.

Particularly in light of the fact that double distributivity appears to be important for getting a good theory of matroids over a Krasner hyperfield [1], it appears to be worthwhile to explore the theory of double distributive hypersemirings and their polynomials.

*Connections to topos theory*

One of the most elegant definitions of a topos is: A category with finite limits and power objects. In fact, every topos has an associated power monad, whose algebras are the internal suplattices. Since the power monad appears to be so fundamental in topos theory, it strikes us as inevitable that there should be a dialogue between the present work and the broad theory of toposes.

As a concrete example, the definition of quantic semirings appears to make sense in an arbitrary topos, and the construction of the universal quantale $\mathcal{M}(\mathcal{S})$ is most likely valid in any topos with a natural numbers object.

On a more speculative note, there has been a lot of interesting work done on the category of relations in a topos (e.g. [5]). In the category of sets, and certainly in greater generality, a relation on $A \times B$ can be identified with a hyper-map $A \to \mathcal{P}(B)$. It would appear fruitful to make a precise comparison between hyperalgebra and the algebra of relations in this context.